\documentclass[12pt]{amsart}
\usepackage{amsthm}
\usepackage{amstext}
\usepackage{amssymb}
\usepackage{mathrsfs}
\usepackage{nicematrix}

\usepackage[english]{babel}
\usepackage{tikz-cd}
\usepackage{enumerate}
\usepackage{xfrac}
\usepackage{mathtools}
\usepackage{xcolor}
\usepackage[all]{xy}
\usepackage{hyperref}
\usepackage{orcidlink}

\usepackage[T1]{fontenc}

\usepackage{soul}
\usepackage{todonotes}

\theoremstyle{plain}
\newtheorem{theorem}{Theorem}[section]
\newtheorem{proposition}[theorem]{Proposition}
\newtheorem{corollary}[theorem]{Corollary}
\newtheorem{lemma}[theorem]{Lemma}

\theoremstyle{definition}
\newtheorem{definition}[theorem]{Definition}

\newtheorem{example}[theorem]{Example}
\newtheorem{remark}[theorem]{Remark}

\theoremstyle{remark}

\numberwithin{equation}{section}

\newcommand{\N}{\mathbb N}

\newcommand{\R}{\mathbb R}
\newcommand{\C}{\mathbb C}
\newcommand{\K}{\mathbb K}

\newcommand{\Mos}{\mathsf{Mos}}
\newcommand{\w}{\mathsf{w}}
\newcommand{\s}{\mathsf{s}}
\newcommand{\srs}{\mathsf{srs}}
\newcommand{\m}{\mathsf{m}}
\newcommand{\ms}{\mathsf{ms}}

\def\rm{\mathrm}

\newcommand{\ICC}{\mathscr{C}}
\newcommand{\loc}{\mathrm{loc}}

\newcommand{\ILL}{\mathscr{L}}

\newcommand{\IMM}{\mathscr{M}}

\newcommand{\IHH}{\mathscr{H}}
\newcommand{\IKK}{\mathscr{K}}
\newcommand{\INN}{\mathscr{N}}
\newcommand{\IVV}{\mathscr{V}}
\newcommand{\IWW}{\mathscr{W}}
\newcommand{\IDD}{\mathscr{D}}
\newcommand{\IEE}{\mathscr{E}}
\newcommand{\IFF}{\mathscr{F}}

\newcommand{\IGG}{\mathscr{G}}

\newcommand{\KK}{\mathscr{K}}

\newcommand{\ISS}{\mathscr{S}}

\newcommand{\ICl}{\mathsf{Cl}}
\newcommand{\Isaplus}{\mathsf{S}^+}
\newcommand{\Isad}{\mathsf{S}}

\newcommand{\set}[2]{\left\{#1\mid #2\right\}}
\newcommand{\elles}{\ILL^{\s}}
\newcommand{\elless}[2]{\elles\left(\mathscr{#1},\mathscr{#2}\right)}
\newcommand{\elle}[2]{\ILL\left(\mathscr{#1},\mathscr{#2}\right)}
\newcommand{\form}[1]{\mathfrak{#1}}
\newcommand{\net}[3]{( #1_{#2};\,#2\in #3)}
\newcommand{\snet}[3]{( #1_{#2})_{#2\in #3}}
\title[]{Strongly continuous fields of operators over varying Hilbert spaces}

\author[A.~BenAmor]{Ali BenAmor$^1$}
\address{$^1$High school for transport and logistics, University of Sousse, Sousse, Tunisia}
\email{ali.benamor@ipeit.rnu.tn}

\author[B.~Güneysu]{Batu Güneysu$^2$\orcidlink{0009-0005-9541-4074}}
\address{$^2$Chemnitz University of Technology, Faculty of Mathematics, 09107 Chemnitz, Germany}
\email{batu.gueneysu@math.tu-chemnitz.de}

\author[T.~Kalmes]{Thomas Kalmes$^2$\orcidlink{0000-0001-7542-1334}}
\email{thomas.kalmes@math.tu-chemnitz.de}

\author[P.~Stollmann]{Peter Stollmann$^2$}
\email{peter.stollmann@math.tu-chemnitz.de}

\begin{document}

\setlength{\parindent}{0pt}

\begin{abstract} After introducing a natural notion of continuous fields of locally convex spaces, we establish a new theory of strongly continuous families of possibly unbounded self-adjoint operators over varying Hilbert spaces. This setting allows to treat operator families defined on bundles of Hilbert spaces that are not locally trivial (such as e.g.~the tangent bundle of Wasserstein space), without referring to identification operators at all.
\end{abstract}

\maketitle

\section{Introduction}
In this paper we introduce a general framework that allows to study continuity - in the strong operator topology - for families of operators over varying spaces. As such, this includes convergence of operators that live on different spaces, a topic of fundamental importance in various fields of analysis, e.g., in numerical analysis, partial differential equations and in stochastic analysis.

With the work of Gromov, Cheeger/Colding and others on limits of sequences of Riemannian manifolds, there has been a renewed interest, in particular in spectral consequences, e.g., for sequences of Laplacians over varying manifolds, making it clear that it is necessary to obtain a set-up that allows to pass from bounded operators to unbounded self-adjoint operators. We mention the fundamental work of Kuwae \&  Shioya, \cite{kuwae}, that introduces the concept of convergence of spectral structures for nets over varying spaces allowing to treat the aforementioned cases.

Another source of applications is the study of energy forms on certain fractals or more general Dirichlet forms by discrete approximation, see \cite{kigami,alonso,hinz} and the literature cited there, as well as the study of continuity aspects for Dirichlet forms and Krein-Feller operators on varying Lebesgue spaces \cite{freiberg-minorics,ehnes-hambly,benamor-kjm}. We also mention the relevance of the framework for approximation of semi-groups of operators as well as their generators  \cite{kato,seidman,ethier-kurz}.

It is rather clear from the beginning that a coherent notion of continuity or convergence for families of vectors from varying spaces lies at the heart of the matter. Once such a notion is established, a family of bounded operators should be viewed as continuous, if it maps continuous sections of vectors to continuous sections of vectors. In the Hilbert space context one can then use functional calculus to pass from unbounded self-adjoint operators to bounded ones.

Apart from Stummel's approach in \cite{stummel1}, who starts with an axiomatic framework, in the existing literature, continuous families of vectors are defined through identification operators between the different spaces. These can be thought of as providing a concept of constant families and continuous families are then naturally identified as those that are close to constant ones.

The starting point of our work is an axiomatic framework by Godement \cite{god} from the 1950's that had been extended by Dixmier \& Douady, \cite{dix}, in the early 1960's. These authors introduce the concept of a continuous field of Banach, resp.~Hilbert spaces over a topological space. It is based on the insight of what properties the set of continuous sections has to share in order to allow a concise theory. At the same time it brings into the game the very useful language of fiber bundles, but we are quick to add that there is not and should not be a requirement of local triviality involved.

Let us point out the main new contributions of this paper:

In Section 2, we extend the notion of continuous fields to  locally convex spaces. This allows to introduce the weak topology on the total space of a continuous field of Hilbert spaces. We can show that many fundamental properties of the interplay between the strong and the weak topology on a single Hilbert space carry over to total spaces, see Section \ref{sec:weakandstrong}. This, in turn, is important in the study of self-adjoint operators and their associated forms in Section \ref{sec:self-adjointops}.  We point out that for the typical arguments one uses, e.g.~passing to subnets, it is essential to actually work with a topology. For a general concept of net convergence this is not necessarily possible. Therefore, the generalization to locally convex spaces is crucial for the later analysis and a major progress.  

In Section \ref{sec:opsbetweencontfields}, we study 
the field of bounded operators between locally convex fields. Provided the underlying space is locally compact, we are able to prove that the natural ad hoc notion of strong continuity mentioned earlier gives the structure of a continuous field of locally convex spaces on the total space of bounded operators between two continuous fields of Hilbert spaces.

For fields of Banach spaces, we introduce in \ref{sec:closedops} the notion of G-continuity of sections, resp.~G-convergence of nets on the total space of closed operators. It provides new insight even for the case of operators between two fixed Banach spaces: the framework allows to pass from  unbounded closed operators to the induced bounded ones by passing to the graph norm. Even if the family we started with lives on a fixed space, we automatically arrive at a family of Banach spaces, in general. G-convergence is compatible with passing to the family of inverse operators (provided they exist). Thus, when specialized to families of self-adjoint operators over a continuous field of Hilbert spaces, we recover the familiar notion of strong resolvent convergence, as seen in Section \ref{sec:self-adjointops}.

In this latter section, we study nonnegative self-adjoint operators through their associated closed forms: the family of form domains endowed with the form norms provides another natural candidate for a continuous field of Hilbert spaces and we are able to show that this is indeed the case, if and only if the operator domains induce a continuous field in the sense of Section \ref{sec:closedops}, and in turn this is equivalent to strong resolvent convergence. This provides a very natural framework for and new insights into Mosco convergence, even for the classical case of forms on one Hilbert space. 

To sum up, we establish total spaces of a continuous fields of Hilbert spaces as a natural framework for perturbation theory. In contrast to earlier approaches it is therefore not necessary to produce mappings between different spaces to study convergence in one fixed Hilbert space, rather one can work in the total space of a continuous field of Hilbert spaces with pretty much the same machinery at hand as in a Hilbert space.   

In the final section, we show how previously defined notions of strong convergence such as e.g.~the mentioned work of Kuwae \& Shioya fit into the context we present here. This gives an abundance of applications of our results. In addition, we give several important examples of continuous fields of locally convex spaces which are not locally trivial in that section, e.g.~the bundle induced by Radon measures on a locally compact space, or the tangent bundle of Wasserstein space. 

Finally, we present two completely new applications of our theory of continuous fields of unbounded self-adjoint operators: firstly, we show that given a continuous family of Riemannian metrics $g^{(x)}$, $x\in X$, on a fixed noncompact manifold $M$ (such as e.g. a Ricci flow), the family of Schrödinger operators
$$
X\times L^2_\loc(M)_+\ni (x,w)\longmapsto -\Delta_{g^{(x)}}+w
$$
is strongly continuous, where $\Delta_{g^{(x)}}$ denotes the Laplace-Beltrami operator induced by $g^{(x)}$.

Secondly, assume $T$ is a locally compact space, $m$ a Radon measure with full support thereon, and $H$ is the self-adjoint operator associated with a regular Dirichlet form $\form{h}$ on $L^2(T,m)$. Denote with $\IMM^R_0(T)$ the space of positive Radon measures on $X$ which are absolutely continuous with respect to the capacity induced by  $\form{h}$, equipped with the weak-*-topology. We then prove that the assignment 
$$
\IMM^R_0(T)\ni \mu \longmapsto  H+\mu
$$
is strongly continuous at each point of $\IMM^R_0(T)$ which is singular with respect to the reference measure $m$.

While the first example naturally lives on continuous field of Hilbert spaces, the second example lives on one Hilbert space, however, the proof of the asserted continuity relies on our results for continuous fields in a very crucial way (via the bundle of Radon measures on a locally compact space).

We have not regarded generalizations of \emph{norm convergence} in this paper, see \cite{post} for an overview of possible concepts as well as \cite{beckusbellissard}. However, we want to underline that the framework of continuous fields is expected to provide a very natural conceptual basis for studying norm continuous sections of bounded operators.\vspace{3mm}


\emph{The last named author wishes to dedicate this work to the memory of his late parents, Margret and Alfred Stollmann in gratitude for their love and support.}

\section{Continuous fields of locally convex spaces}
\label{sec:fields}
We understand all our linear spaces and function spaces over the fixed field $\K$, where $\K$ denotes either the complex or the real numbers. Following the usual abuse of notation, we denote the scalar product on any Hilbert space by $\left\langle \cdot,\cdot\right\rangle$, and the norm on a normed space with $\left\|\cdot\right\|$. Given a set $X$ and a family $\IVV:=(\IVV_x;x\in X)$ of sets, $\pi:\IVV\to X$ denotes the canonically given projection $\zeta\mapsto x$, if $\zeta\in \IVV_x$. Whenever convenient, we will assume the $\IVV_x$ to be pairwise disjoint, so there is no harm in identifying
$$
\IVV=\bigsqcup_{x\in X}\IVV_x,
$$
where the latter stands for the disjoint union. Concerning families, we will change between index and argument notation $\IVV(x)$ to ease notation. 
\begin{definition} One calls $\IVV_x=\pi^{-1}(\{x\})$ the \emph{fibers of $\IVV$}, and the \emph{space of sections of $\IVV$}, to be denoted with $\Gamma(\IVV)$, is defined to be the set of all maps $f:X\to \IVV$ satisfying $\pi \circ f=\mathrm{id}_X$.
\end{definition}

If all fibers of $\IVV$ are linear spaces, then the space of sections of $\IVV$ becomes a left-module over the ring of functions on $X$ (in particular, a linear space).

In what follows, $(X, \mathfrak{T})$ will be a topological space. In the sequel the set of neighborhoods of a point $x$ in a topological space is denoted by $\mathfrak{U}(x)$, the topology in question being understood.

\begin{definition} \label{cflcs}
\begin{enumerate}
 \item A \emph{field of locally convex spaces over $X$ with index set $I$} is a pair $(\IVV, p)$, where for every $x\in X$ the space $\IVV_x$ is a locally convex space and $p=(p_i(x); i\in I, x\in X)$ is such that $(p_i(x); i\in I)$ is a family of seminorms on $\IVV_x$ that generates its topology. We point out that by the latter we mean that for every $x\in X$ the set $\{B_{p_i(x)}(0,\epsilon); i\in I, \varepsilon>0\}$ of open balls of radius $\epsilon>0$ around zero in $\IVV_x$ with respect to the seminorm $p_i(x)$, $B_{p_i(x)}(0,\epsilon)$, is a zero neighborhood basis in $\IVV_x$. Additionally, we emphasize that we do not require locally convex spaces to be Hausdorff.
 \item A subspace $\ICC\subset \Gamma(\IVV)$ is called \emph{continuous} for $(\IVV, p)$, if for all $f\in\ICC$ and $i\in I$, the function $x \mapsto p_i(x)(f(x))$ is continuous (we will use the shorthand $p_i(\cdot)(f(\cdot))$ or $p_i(f)$ for functions of the latter type in what follows).
 
 \item A subspace $\ICC\subset \Gamma(\IVV)$ is called 
 \emph{locally uniformly closed}, if it satisfies the property: 
    \begin{itemize}
    \item[\textrm{(LUC)}]
 Every $g\in \Gamma(\IVV)$ with
 \begin{align}\label{luc}
 \mbox{for  any }i\in I, x\in X\mbox{ and }\varepsilon>0&, \mbox{ there exist }f \in\ICC\mbox{ and }U\in \mathfrak{U}(x)
 \mbox{ such that}\\
            p_i(y)\left(g(y)-f(y)\right)&\le \varepsilon\mbox{ for all }y\in U,\nonumber
 \end{align}
            satisfies $g\in \ICC$.
            \end{itemize}
 \item  A subspace $\ICC\subset \Gamma(\IVV)$ is called \emph{saturated}, if for all $x\in X$, $\zeta\in \IVV_x$ there exists $f\in \ICC$ with $f(x)=\zeta$.
 \item A \emph{continuous field of locally convex spaces over $X$} is a triple $(\IVV,\ICC,p)$ so that $(\IVV,p)$ is a field of locally convex spaces and $\ICC$ is a continuous subspace that is locally uniformly closed and saturated.
\item We speak of a [\emph{continuous}] \emph{field} $\IFF$ \emph{of Fr\'{e}chet spaces}, if all the $\IFF_x$ share this property and we assume that in this case $I=\N$, $p_n(x)\leq p_{n+1}(x)$ for all $x\in X, n\in \N$. We denote the canonical translation invariant metric on $\IFF_x$ defined via the $(p_n(x); n\in\N)$ and inducing the topology of $\IFF_x$ by $d_{\IFF_x}$, i.e.
\begin{equation}\label{f-d}
d_{\IFF_x}(\zeta,\zeta'):=\sum^\infty_{n=1}\frac{1}{2^n}\frac{p_n(x)(\zeta-\zeta')}{1+p_n(x)(\zeta-\zeta')}\quad (\zeta,\zeta'\in\IFF_x)
\end{equation}
and observe that $d_{\IFF_x}(\lambda \zeta,\lambda \zeta')\leq d_{\IFF_x}(\zeta,\zeta')$ for $\zeta,\zeta'\in\IFF_x$ and $\lambda\in\K, |\lambda|\leq 1$. It is readily seen that  the function $x\mapsto d_{\IFF_x}(f(x), g(x))$ is continuous for sections $f,g$ from a continuous subspace $\ICC$ for $\IFF$. 
\item If all the spaces $\IVV_x$ are normed, we assume that $I$ consists of only one element and omit indices; as mentioned above, we then simply write $\|\cdot\|$ for the norm on $\IVV_x$, suppressing the reference to $x$. 
\item Finally, we speak of a [\emph{continuous}] \emph{field} \emph{of Banach spaces} and \emph{Hilbert spaces}, respectively, if all the fibers share the respective property, with the above convention concerning the norm and the inner product. In these cases, we simply denote the continuous field by $(\IVV, \ICC)$.
\end{enumerate}
\end{definition}

\begin{remark} 
\begin{itemize}
\item[(1)] Definition \ref{cflcs} is modeled on \cite{god,dix} (see also \cite{bos}), where the fibers are Hilbert or Banach spaces. Its generalization to locally convex spaces is considered here mainly to include certain fields of spaces that arise naturally: fields of Hilbert spaces endowed with the weak topology and fields of bounded operators endowed with the strong operator topology.
\item[(2)] Note that the index set $I$ that is used to label the seminorms is the same for all fibers. Of course, this could be achieved for any family of locally convex spaces. However, the choice of $p_i$ is important in terms of continuous subspaces as it amounts to a sort of continuous section of the family of continuous seminorms for the different $\IVV_x$.
\item[(3)] The property (LUC) is very natural in that its analog is clearly satisfied for the space of continuous mappings $C(X;\IVV_0)$ into a single locally convex space.
\item[(4)] Expanding on the last point, continuity of a mapping $f:X\to\IVV_0$ at a point $x_0$ is obviously equivalent to the continuity of $p_i(f-f(x_0))$ at $x_0$ for all $i\in I$. This observation is implicit in the proof of the following Lemma.
\end{itemize}
\end{remark}

\begin{lemma}\label{Lemma2.4}
 Let $(\IVV, p)$ be a field of locally convex spaces with index set $I$ over $X$, and let $\ICC$ be a continuous subspace for $(\IVV,p)$. Then 
 \begin{itemize}
  \item[\textrm{(1)}] 
  $$\overline{\ICC}:=\set{g\in\Gamma(\IVV)}{g\mbox{  satisfies \eqref{luc}}}$$
  defines the smallest locally uniformly closed continuous subspace containing $\ICC$.
  \item[\textrm{(2)}] Assume that $\ICC$ satisfies the property 
    \begin{itemize}
    \item[\textrm{(LUC')}] For all $g\in\Gamma(\IVV)$ satisfying 
$$
p_i(\cdot)(g(\cdot)-f(\cdot))\in C(X)\quad\text{for all $i\in I$, $f\in\ICC$}, 
$$
one has $g\in\ICC$.
\end{itemize}
Then $\ICC$ is locally uniformly closed. Conversely, if $\ICC$ is locally uniformly closed and saturated, then $\ICC$ satisfies \textrm{(LUC')}.
\item[\textrm{(3)}] If $\ICC$ is locally uniformly closed, it is a left $C(X)$ module. 
 \end{itemize}
\end{lemma}

\begin{proof}
(1) It is clear that every locally uniformly closed continuous subspace containing $\ICC$ must also contain $\overline{\ICC}$, so it remains to show that $\overline{\ICC}$ is a locally uniformly closed continuous subspace. Obviously, $\overline{\ICC}$ is a subspace of $\Gamma(\IVV)$ and that it is locally uniformly closed is easy to see: assume that $g\in \Gamma(\IVV)$ satisfies \eqref{luc} for $\overline{\ICC}$ in place of $\ICC$. This means, that for given, $\varepsilon>0$, $i\in I$, there are $f \in\overline{\ICC}$ and $V\in \mathfrak{U}(x)$ so that 
$$p_i(y)\left(g(y)-f(y)\right)\le \varepsilon\mbox{ for all }y\in V.$$
The definition of $\overline{\ICC}$ now gives $U\subset V\in \mathfrak{U}(x)$ and $f'\in\ICC$ so that 
$$p_i(y)\left(f'(y)-f(y)\right)\le \varepsilon\mbox{ for all }y\in U.$$
The triangle inequality gives that  
$$p_i(y)\left(g(y)-f'(y)\right)\le 2\varepsilon\mbox{ for all }y\in U,$$
so that $g\in\overline{\ICC}$.

To show that $\overline{\ICC}$ is a continuous subspace, we first observe that for every $g_1,g_2,g_3\in\Gamma(\IVV)$, each $x,y\in X$, and every $i\in I$ it holds
\begin{eqnarray*}
    p_i(y)\left(g_1(y)-g_3(y)\right)&\leq&p_i(y)\left(g_1(y)-g_2(y)\right)+p_i(y)\left(g_2(y)-g_3(y)\right)-p_i(x)\left(g_2(x)-g_3(x)\right)\\
    &&+p_i(x)\left(g_2(x)-g_3(x)\right)\\
    &\leq&p_i(y)\left(g_1(y)-g_2(y)\right)+\left|p_i(y)\left(g_2(y)-g_3(y)\right)-p_i(x)\left(g_2(x)-g_3(x)\right)\right|\\
    &&+p_i(x)\left(g_2(x)-g_1(x)\right)+p_i(x)\left(g_1(x)-g_3(x)\right),
\end{eqnarray*}
implying
\begin{eqnarray*}
    p_i(y)\left(g_1(y)-g_3(y)\right)-p_i(x)\left(g_1(x)-g_3(x)\right)&\leq&p_i(y)\left(g_1(y)-g_2(y)\right)+p_i(x)\left(g_2(x)-g_1(x)\right)\\
    &&+\left|p_i(y)\left(g_2(y)-g_3(y)\right)-p_i(x)\left(g_2(x)-g_3(x)\right)\right|.
\end{eqnarray*}
Combining the previous inequality with itself when the roles of $x$ and $y$ are interchanged, we obtain
\begin{eqnarray}\label{eq:auxiliary equation 1 for equivalence (C3) and (C3')}
    \left|p_i(y)\left(g_1(y)-g_3(y)\right)-p_i(x)\left(g_1(x)-g_3(x)\right)\right|&\leq&p_i(y)\left(g_1(y)-g_2(y)\right)+p_i(x)\left(g_2(x)-g_1(x)\right)\nonumber \\
    &&+\left|p_i(y)\left(g_2(y)-g_3(y)\right)-p_i(x)\left(g_2(x)-g_3(x)\right)\right|
\end{eqnarray}
for every $g_1,g_2,g_3\in\Gamma(\IVV)$, each $x,y\in X$ and every $i\in I$. 

Now, let $g\in\overline{\ICC}$ be arbitrary. Fixing $i\in I$ and $x\in X$ it follows from \eqref{eq:auxiliary equation 1 for equivalence (C3) and (C3')} evaluated for $g_1=g, g_2=f$ and $g_3=0$ that
\begin{eqnarray*}
    |p_i(y)\left(g(y)\right)-p_i(x)\left(g(x)\right)|&\leq& p_i(y)\left(g(y)-f(y)\right)+p_i(x)\left(f(x)-g(x)\right)\\
    &&+|p_i(y)\left(f(y)\right)-p_i(x)\left(f(x)\right)|
\end{eqnarray*}
holds for every $f\in\ICC$ and $y\in X$. The continuity of $y\mapsto p_i(y)\left(g(y)\right))$ in $x$ follows from this inequality since $g$ satisfies \eqref{luc}. Hence, $\overline{\ICC}$ is a continuous subspace for $(\IVV,p)$.

(2)  We first prove that (LUC') implies (LUC). Hence, fix $g\in \Gamma(\IVV)$ as in \eqref{luc}. By (LUC'), in order to verify $g\in\mathscr{C}$, it suffices to show the continuity of $p_i(\cdot)\left(g(\cdot)-f(\cdot)\right)$ for every choice of $f\in \mathscr{C}$ and all $i\in I$. Thus, fix $f\in\mathscr{C}, i\in I$ and let $x\in X, \varepsilon>0$ be arbitrary. By \eqref{luc}, there are $f'\in\mathscr{C}$ and $U\in \mathfrak{U}(x)$ such that
\begin{equation}\label{eq:auxiliary equation 2 for equivalence (C3) and (C3')}
    p_i(y)\left(g(y)-f'(y)\right)\leq\varepsilon,\quad (y\in U).
\end{equation}
Additionally, $f-f'\in\mathscr{C}$, so there is $V\in\mathfrak{U}(x)$ with
\begin{equation}\label{eq:auxiliary equation 3 for equivalence (C3) and (C3')}
    \left|p_i(y)\left(f(y)-f'(y)\right)-p_i(x)\left(f(x)-f'(x)\right)\right|\leq\varepsilon,\quad (y\in V).
\end{equation}
Combining \eqref{eq:auxiliary equation 1 for equivalence (C3) and (C3')} for $g_1=g, g_2=f, g_3=f'$, \eqref{eq:auxiliary equation 2 for equivalence (C3) and (C3')}, and \eqref{eq:auxiliary equation 3 for equivalence (C3) and (C3')}, it follows
\[
     \left|p_i(y)\left(g(y)-f'(y)\right)-p_i(x)\left(g(x)-f'(x)\right)\right|\leq 3\varepsilon,\quad (y\in U\cap V)
\]
which proves the continuity of $p_i(\cdot)\left(g(\cdot)-f(\cdot)\right)$ in $x$. 

In order to show that (LUC) also implies (LUC') for saturated continuous subspaces, let $g\in\Gamma(\IVV)$ such that $p_i(\cdot)\left(g(\cdot)-f(\cdot)\right)$ is continuous for every $i\in I, f\in\mathscr{C}$. To conclude $g\in\mathscr{C}$, by (LUC), it suffices to show that for any $i\in I, x\in X$ and $\varepsilon>0$, there are $f'\in\mathscr{C}$ and $U\in \mathfrak{U}(x)$ such that $p_i(y)\left(g(y)-f'(y)\right)\leq \varepsilon$ for $y\in U$. Fix $i\in I, x\in X$ and $\varepsilon>0$. Since $\mathscr{C}$ is saturated, there is $f'\in\mathscr{C}$ with $f'(x)=g(x)$. The continuity of $p_i(\cdot)\left(g(\cdot)-f'(\cdot)\right)$ in $x$ grants the existence of $U\in\mathfrak{U}(x)$ with
\[
    \left|p_i(y)\left(g(y)-f'(y)\right)-p_i(x)\left(g(x)-f'(x)\right)\right|\leq \varepsilon,\quad (y\in U).
\]
The desired inequality $p_i(y)\left(g(y)-f'(y)\right)\leq \varepsilon$ for $y\in U$ follows herefrom with the aid of \eqref{eq:auxiliary equation 1 for equivalence (C3) and (C3')} for $g_1=g_2=g$ and $g_3=f'$. 

(3) is again easy: for $f\in\ICC$ and $\varphi\in C(X)$ we can approximate $\varphi\cdot f$ by $\varphi(x)\cdot f$
in a small enough neighborhood of $x\in X$, and $\varphi(x)\cdot f\in\ICC$ by linearity.
\end{proof}

Next, we extend a useful result to the locally convex setting that provides a sufficient condition, that a space of sections generates a strongly continuous field 'by closure'.

\begin{lemma}\label{completion}
  Assume that $(\IFF,p)$ is a field of Fr\'{e}chet spaces over $X$ and $\ICC$ is a continuous subspace for $(\IFF,p)$. If $\set{f(x)}{f\in\ICC}$ is dense in $\IFF_x$ for any $x\in X$, then $\overline{\ICC}$ from Lemma \ref{Lemma2.4} (1) is saturated. In particular, $(\IFF, \overline{\ICC},p)$ is a continuous field of Fr\'{e}chet spaces over $X$.
\end{lemma}
\begin{proof}
By Lemma \ref{Lemma2.4} (1), $\overline{\ICC}$ is locally uniformly closed, which explains the \emph{in particular} part. It remains to show that $\overline{\ICC}$ is saturated. Fix $x\in X$. It suffices to show that 
$$
\overline{\ICC}_x:=\set{f(x)}{f\in\overline{\ICC}}\subset \IFF_x
$$
is closed. Thus, let $\zeta$ be in the closure of $\overline{\ICC}_x$ in $\IFF_x$. Since the latter is a Fr\'echet space, there is a sequence $(\zeta_l)_{l\in\N}$ in $\overline{\ICC}_x$ converging to $\zeta$. For $n\in\N$ let $l_n$ be such that
$$p_n(x)\left(\zeta-\zeta_l\right)<2^{-(n+2)}\text{ for all }l\geq l_n.$$
Without loss of generality we assume $l_{n+1}>l_n$, $n\in\N$. Because $\overline{\ICC}$ is a subspace of $\Gamma(\IFF)$, there are $g, f_m\in\overline{\ICC}$, $m\in \N$, such that $g(x)=\zeta_{l_1}$, $f_m(x)=\zeta_{l_{m+1}}-\zeta_{l_m}$. For each $m\in\N$ we set
$$
f'_{m}(\cdot):=\left(\max_{j=1,\dots,m}\left(1,\frac{p_j(\cdot)(f_m(\cdot))}{p_j(x)(\zeta_{l_{m+1}}-\zeta_{l_m})}\right)\right)^{-1}f_m(\cdot)\in \overline{\ICC},
$$
so that $f'_m(x)=\zeta_{l_{m+1}}-\zeta_{l_m}$. By our general assumption on fields of Fr\'echet spaces, for every $y\in X$ the sequence of seminorms $(p_n(y))_n$ is increasing so that for every $n,j,k\in \N$ with $j\geq n$
\begin{eqnarray}\label{eq:closure 1}
    p_n(y)\left(\sum_{m=j}^{j+k} f'_m(y)\right)&\leq& \sum_{m=j}^{j+k} p_n(y)\left(f'_m(y)\right)\leq \sum_{m=j}^{j+k} p_n(x)\left(\zeta_{l_{m+1}}-\zeta_{l_m}\right)\nonumber\\
    &\leq& \sum_{m=j}^{j+k} p_m(x)\left(\zeta_{l_{m+1}}-\zeta_{l_m}\right)\leq 2^{-j}.
\end{eqnarray}
It follows that $\sum_{m=1}^\infty f'_m(y)$ converges in $\IFF_y$ for every $y\in X$. Thus, $h:=g+\sum_{m=1}^\infty f'_m\in\Gamma(\IFF)$ and $h(x)=g(x)+\sum_{m=1}^\infty f'_m(x)=\zeta$.

In order to show that $\zeta\in \overline{\ICC}_x$, which will entail that $\overline{\ICC}_x$ is closed in $\IFF_x$, we shall prove $h\in\overline{\ICC}$. For this, we fix $n\in\N, y\in X$ and $\epsilon>0$ and we choose $j\in\N$, $j\geq n,$ such that $2^{-j}<\epsilon/2$. On the one hand, by \eqref{eq:closure 1}
$$p_n(y)\left(h(y)-g(y)-\sum_{m=1}^{j} f'_m(y)\right)<\epsilon/2, \quad y\in X,$$
on the other hand, as $g+\sum_{m=1}^j f'_m\in \overline{\ICC}$, there are $U\in \mathfrak{U}(y)$ and $f\in \ICC$ satisfying
$$p_n(z)\left(g(z)+\sum_{m=1}^j f'_m(z)-f(z)\right)<\epsilon/2, \quad z\in U,$$
from which we obtain $h\in \overline{\ICC}$ as desired.
\end{proof}

%



Given a continuous field of locally convex spaces $(\IVV,\ICC,p)$ with index set $I$ over $X$, we consider the sets

$$
U(V,f,i,\epsilon):=\set{\zeta\in \IVV}{\pi(\zeta)\in V, p_i(\pi(\zeta))\big(\zeta-f(\pi(\zeta))\big)<\epsilon}\subset \IVV,
$$
where $V\subset X$ is open, $f\in \ICC$, $i\in I$, and $\epsilon>0$. 
    
In case of a continuous field $(\IVV,\ICC)$ of Banach spaces over $X$, as then $(p_i(x))_{i\in I}$ consists of a single norm for each $x\in X$, we simply write
$$
U(V,f,\epsilon):=\set{\zeta\in \IVV}{\pi(\zeta)\in V, \|\zeta-f(\pi(\zeta))\|<\epsilon}.
$$

\begin{proposition}\label{prop:properties of the topology}
Assume $(\IVV,\ICC,p)$ is a continuous field of locally convex spaces with index set $I$ over $X$. 
\begin{itemize}
    \item [(1)] The sets $U(V,f,i,\epsilon)$ above form the base of a topology on $\IVV$ (and $\IVV$ will always be equipped with this topology in the sequel). This topology is Hausdorff if $X$ as well as all the fibers $\IVV_x$ are so.
    \item [(2)] The projection $\pi:\IVV\to X$ is continuous and open, and the fibers carry the subspace topology. Additionally, $p_i:\IVV\rightarrow \R, \zeta\mapsto p_i(\pi(\zeta))(\zeta)$ is continuous for every $i\in I$.
    \item [(3)] A section $f\in\Gamma(\IVV)$ is in $\ICC$, if and only if $f$ is a continuous map. In other words, 
    $$
    \ICC=\Gamma(\IVV)\cap C(X,\IVV).
    $$
    \item[(4)] Given $U\subset X$ open, the subspace
    $$
    \IVV|_U:=\bigsqcup_{x\in U} \IVV_x\subset \IVV
    $$
    together with 
    $$
    \ICC_U:=\Gamma(\IVV|_U)\cap C(U,\IVV)
    $$
    is a continuous field of locally convex (resp. Banach, resp. Hilbert) spaces.
    \item[(5)] The assignment $U\mapsto \ICC_U$ is a sheaf of modules with respect to the sheaf of continuous functions on $X$.
    \item[(6)] Let $\zeta\in U(V,f,i,\epsilon)$ and assume that $g\in\ICC$ satisfies
    $$
    p_i(\pi(\zeta))\left( g(\pi(\zeta))-f(\pi(\zeta))\right)<\epsilon .
    $$
    Then there is $\hat{V}\in\mathfrak{U}(\pi(\zeta))$ and $\hat{\epsilon}>0$ with
    $$
    U(\hat{V},g,i,\hat{\epsilon})\subset U(V,f,i,\epsilon) .
    $$
\end{itemize}
\label{strong-top}
\end{proposition}

\begin{proof} (1): (a) Let $\zeta^0\in\IVV$, $x=\pi(\zeta^0)$. Pick $f\in\ICC$ such that $f(x)=\zeta^0$. It follows that for every open subset $V\subset X$ containing $x$ and every $\epsilon>0$, $i\in I$ we have $\zeta^0\in U(V,f,i,\epsilon)$.\\
(b) Let 
$$
\zeta^0\in U(V_1,f_1,i_1,\epsilon_1) \cap U(V_2,f_2,i_2,\epsilon_2)
$$ 
and set $x:=\pi(\zeta^0)$. It follows that $x\in V_1\cap V_2$ as well as
\[
p_{i_1}(x)\left( \zeta^0 - f_1(x) \right ) <\epsilon_1,\ p_{i_2}(x)\left( \zeta^0 - f_2(x) \right ) <\epsilon_2.
\]
Because $(p_i(x); i\in I)$ generates the topology of $\IVV_x$, there are $i_3\in I$ and $\epsilon_3>0$ such that
\[
B_{p_{i_3}(x)}(0,\epsilon_3)\subset B_{p_{i_1}(x)}(0,\epsilon_1)\cap B_{p_{i_2}(x)}(0,\epsilon_2)
\]
which easily implies
\[
p_{i_j}(x)\leq \frac{\epsilon_j}{\epsilon_3} p_{i_3}(x), \quad j=1,2.
\]
In particular, there is $C>0$ with $p_{i_j}(x)\leq C p_{i_3}(x)$. Pick $\epsilon\in(0,\min\{\epsilon_1,\epsilon_2\})$ so small that 
\[
p_{i_j}(x)\left( \zeta^0 - f_j(x) \right ) + \epsilon <\epsilon_j, \quad j=1,2,
\]
and choose $f\in\ICC$ such that $f(x) = \zeta^0$. By continuity of $p_{i_3}(\cdot)\left( f(\cdot) - f_j(\cdot) \right )$, $j=1,2$, we find $V\subset V_1\cap V_2, V$ open, containing $x$, such that 
\[
p_{i_j}(y)\left( f(y) - f_j(y) \right )< (\min\{\epsilon_1,\epsilon_2\} - \epsilon)/C, \quad y\in V, j=1,2. 
\]
Obviously, $\zeta^0\in U(V,f,i_3,\epsilon/C)$. Let us show that 
$$
U(V,f,i_3,\epsilon/C)\subset U(V_1,f_1,i_1,\epsilon_1) \cap U(V_2,f_2,i_2,\epsilon_2).
$$
Indeed, for $\zeta\in U(V,f,i_3,\epsilon/C)$ it holds $y:=\pi(\zeta)\in V\subset V_1\cap V_2$. Therefore, by the choice of $V$, 
\begin{align*}
p_{i_j}(y)\left( \zeta- f_j(y) \right ) &\leq p_{i_j}(y)\left( \zeta- f(y) \right ) + p_{i_j}(y)\left( f(y) - f_j(y) \right ) \\
& \leq  C p_{i_3}(y)\left( \zeta- f(y) \right ) +  p_{i_j}(y)\left( f(y) - f_j(y) \right ) < \epsilon_j.
\end{align*} 

(c) Let $X$ and each fiber $\IVV_x$ be Hausdorff as well as $\xi_1,\xi_2\in \IVV$, $\xi_1\neq \xi_2$. If $\pi(\xi_1)\neq\pi(\xi_2)$ there are $V_j\in\mathfrak{U}(\xi_j)$, $j=1,2$, with $V_1\cap V_2=\emptyset$. For $f_j\in\ICC$, $j=1,2$, with $f_j(\pi(\xi_j))=\xi_j$ we have $\xi_j\in U(V_j,f_j,i,1)$ and $U(V_1,f_1,i,1)\cap U(V_2,f_2,i,1)=\emptyset$, where $i\in I$ is arbitrary. If $\pi(\xi_1)=\pi(\xi_2)=:x$ there is $i\in I$ with $2\epsilon:=p_i(x)(\xi_1-\xi_2)>0$. For $f_j\in\ICC$, $j=1,2$, with $f_j(\pi(\xi_j))=\xi_j$ we have $\xi_j\in U(V_j,f_j,i,\epsilon)$ and $U(V_1,f_1,i,\epsilon)\cap U(V_2,f_2,i,\epsilon)=\emptyset$.

(2) We have $\pi\left(U(V,f,i,\epsilon)\right)= V$ and hence $\pi$ is an open mapping and the latter identity implies that $\pi$ is continuous.\\
For $x\in X$, $V\subset X$ containing $x$, $f\in \ICC$ and $\zeta = f(x)$, we have for every $i\in I$
\begin{align*}
\IVV_x \cap U(V,f,i,\epsilon) &= \set{\xi\in \IVV_x}{p_i(\pi(\xi)) \left(\xi - f(\pi(\xi)\right)<\epsilon}\\
& = \set{ \xi\in \IVV_x}{p_i(x) \left(\xi - \zeta\right)<\epsilon}
\end{align*}
Consequently these sets form a base of neighborhoods of $\zeta$, by assumptions on $\left( p_i(x)\right)_{i\in I}$. Since any $\zeta\in \IVV_x$ can be attained by a suitable $f\in\ICC$, the subspace topology of each fiber is equal to its original topology.

Now, we fix $i\in I$. For $\zeta\in\IVV$ there is $f\in\IVV$ with $f(\pi(\zeta))=\zeta$. Thus, for $\epsilon>0$ it holds $p_i\left(U(X,f,i,\epsilon\right)\subset B_{|\cdot|}(p_i(\pi(\zeta)(\zeta),\epsilon))$ so that $p_i$ is continuous in $\zeta$.\\
(3) Let $f\in \ICC, x\in X$ and $U(V,g,i,\epsilon)$ a base set containing $\zeta=f(x)$. 
By continuity of $p_i(\cdot) \left( f(\cdot) - g(\cdot) \right)$ we find a neighborhood $W$ of $x$ such that $W\subset V$ and
\[
p_i(y)  \left( f(y) - g(y) \right) <\epsilon \quad\text{for all $y\in W$}.
\]
It follows that $f(y)\in U(V,g,F,\epsilon)$ for all $y\in W$, showing that $f:X\to\IVV$ is continuous.\\
Conversely, let $g\in\Gamma(\IVV)$ be continuous. In particular, for every $i\in I$, $x\in X$, and $\epsilon>0$, there is a neighborhood $V$ of $x$ such that
\[
g(y)\in U(V,f,i,\epsilon)\ \text{for all}\ y\in V. 
\]
Since this is exactly the condition \eqref{luc}, $g\in\ICC$ follows because $\ICC$ is locally uniformly closed. \\
(4) Clearly, with the obvious interpretation, $\IVV|_U$ is a field of locally convex spaces over $U$ with index set $I$ which is a field of Banach spaces, resp.~Hilbert spaces, whenever $\IVV$ is. Applying (2) we see that $\ICC|_U$ is a continuous subspace for $\IVV|_U$. Because $\ICC|_U\subset \ICC_U$, it is saturated.

In order to show that $\ICC_U$ is locally uniformly closed, let $g\in \Gamma(\IVV|_U)$ be such that \eqref{luc} holds with $X$ and $\ICC$ replaced by $U$ and $\ICC_U$, respectively. We have to show $g\in C(U,\IVV)$. Thus, fix $x\in U$, and let $W$ be a neighborhood of $g(x)$ in $\IVV|_U$, i.e.~for a suitable neighborhood $V_1$ of $x$ in $U$, $f_1\in \ICC|_U$, $i\in I$, and $\epsilon>0$ it holds $g(x)\in U(V_1,f_1,i,\epsilon)\subset W$. We choose $\delta>0$ with $p_i(x)(g(x)-f_1(x))+\delta<\epsilon$. Due to \eqref{luc}, for this $i\in I$ and $\delta/3$ in place of $\epsilon$, there are $f_2\in \ICC_U$ and a neighborhood $V_2$ of $x$ in $U$ such that
\[
p_i(y)\left(g(y)-f_2(y)\right)\leq \delta/3,\quad y\in V_2.
\]
Because $f_2-f_1\in \ICC_U=\Gamma(\IVV|_U)\cap C(U,\IVV)$, by (2) for $\IVV|_U$, there is a neighborhood $V_3$ of $x$ in $U$ such that
\[
p_i(y)(f_2(y)-f_1(y))-p_i(x)(f_2(x)-f_1(x))<\delta/3,\quad y\in V_3.
\]
Then, for $y\in V_1\cap V_2\cap V_3$ it holds
\begin{eqnarray*}
    p_i(y)(g(y)-f_1(y))&\leq& p_i(y)(g(y)-f_2(y))+p_i(y)(f_2(y)-f_1(y))\\
    &<&2\delta/3+p_i(x)(f_2(x)-f_1(x))\\
    &\leq& 2\delta/3+p_i(x)(g(x)-f_2(x))+p_i(x)(g(x)-f_1(x))\\
    &\leq& \delta + p_i(x)(g(x)-f_1(x))<\epsilon,
\end{eqnarray*}
where the last inequality holds by the choice of $\delta$. Therefore, $g(y)\in U(V_1,f_1,i,\epsilon)$ for each $y\in V_1\cap V_2\cap V_3$ so that $g$ is continuous in $x$. Thus, $\ICC_U$ is indeed locally uniformly closed which completes the proof of (4).

(5) is obvious.

(6) follows from the continuity of $p_i(g-f)$ and the triangle inequality.
\end{proof}

In particular in applications it may be useful to work with sequences to settle certain topological aspects. Therefore, the following observation is worthwhile to note:

\begin{corollary}
    \label{cor:countable}
    Let $(\IFF,\ICC,p)$ be a continuous field of Fr\'echet spaces over $X$. 
    \begin{itemize}
        \item[(1)] If $X$ is first countable, then so is $(\IFF, \mathfrak{T})$.
        \item[(2)] If $X$ is second countable and $\IFF_x$ is separable for all $x\in X$, then $(\IFF, \mathfrak{T})$ is second countable.
    \end{itemize}    
\end{corollary}
\begin{proof}
    Both assertions follow immediately from Proposition \ref{prop:properties of the topology}, in particular its part (6).
\end{proof}

\begin{remark}\label{rem:natural topology is an initial topology}
    For a continuous field  of locally convex spaces $(\IVV,\ICC,p)$ it is easy to see that the mapping $\K\times\IVV\rightarrow \IVV, (\lambda,\xi)\mapsto \lambda\xi$ is correctly defined and continuous. Similar arguments as the ones used to prove Proposition \ref{prop:properties of the topology} (3) yield the continuity of the mapping $\Delta\rightarrow\IVV, (\xi_1,\xi_2)\mapsto \xi_1+\xi_2$, where
    $$\Delta:=\set{(\xi_1,\xi_2)\in\IVV}{\pi(\xi_1)=\pi(\xi_2)}$$
    is equipped with the (subset topology of the) product topology. Indeed, fix $(\xi_1^0,\xi_2^0)\in \Delta$ and let $x:=\pi(\xi_1^0)=\pi(\xi_2^0)=\pi(\xi_1^0+\xi_2^0)$. For $U(V,f,i,\epsilon)\in\mathfrak{U}(\xi_1^0+\xi_2^0)$ we choose $\delta>0$ with $p_i(x)(\xi_1^0+\xi_2^0-f(x))+\delta<\epsilon$, as well as $f_j\in\ICC$ such that $f_j(x)=\xi_j^0$, $j=1,2$. From $f_1+f_2-f\in\ICC$ we conclude the existence of $V_1\in\mathfrak{U}(x)$ with $V_1\subset V$ and
    $$p_i(y)\left(f_1(y)+f_2(y)-f(y)\right)-p_i(x)\left(f_1(x)+f_2(x)-f(x)\right)<\delta/4,\quad y\in V_1.$$
    Then, for $(\xi_1,\xi_2)\in U(V_1, f_1,i,\delta/4)\times U(V_1, f_2,i,\delta/4)\cap\Delta\in \mathfrak{U}((\xi_1^0,\xi_2^0))$ we derive - keeping in mind that $\pi(\xi_1)=\pi(\xi_2)=\pi(\xi_1+\xi_2)\in V_1$ and $f_j(x)=\xi_j^0$ - 
    \begin{eqnarray*}
        p_i(\pi(\xi_1+\xi_2))(\xi_1+\xi_2-f(\pi(\xi_1+\xi_2)))&\leq&p_i(\pi(\xi_1))(\xi_1-f_1(\pi(\xi_1))+p_i(\pi(\xi_2))(\xi_2-f_2(\pi(\xi_2)))\\
        &&+p_i(\pi(\xi_1))(f_1(\pi(\xi_1))+f_2(\pi(\xi_2))-f(\pi(\xi_1)))\\
        &<&\delta + p_i(x)(f_1(x)+f_2(x)-f(x))<\epsilon,
    \end{eqnarray*}
    proving $\xi_1+\xi_2\in U(V,f,i,\epsilon)$ which shows the claimed continuity.

    As a consequence of the above observation and Proposition \ref{prop:properties of the topology} it follows that for every $i\in I$ and $f\in\ICC$ the function
    $$F(i,f):\IVV\rightarrow\R, \xi\mapsto p_i(\pi(\xi))(\xi-f(\pi(\xi)))$$
    is continuous. From this it follows immediately that the natural topology on $\IVV$ introduced in Proposition \ref{prop:properties of the topology} is the coarsest topology $\mathfrak{T}$ on $\IVV$ for which every function $F(i,f):(\IVV,\mathfrak{T})\rightarrow \R$, $i\in I$, $f\in\ICC$ together with $\pi:(\IVV,\mathfrak{T})\rightarrow X$ are continuous. An analogous description of $\mathfrak{T}$ as an initial topology was given in \cite{grothauswittmann} for the case of sequences of spaces.
\end{remark}

\begin{proposition} \label{ON-Frames}
\emph{(1)} Let $(\IFF,\ICC, p)$ be a continuous field of Fr\'echet spaces over $X$, let $x_0 \in X$ and let $\xi_1, ...,\xi_n\in \IFF_{x_0}$ be linearly independent. Then there is $U\in\mathfrak{U}(x_0)$ and $f_1, ... ,f_n\in\ICC(U)$ such that $f_1(x), ... , f_n(x)$ are linearly independent for all $x\in U$ and $f_j(x_0)=\xi_j$ for all $j=1, ... ,n$.\vspace{1mm}

\emph{(2)} Let $(\IHH,\ICC)$ be a continuous field of Hilbert spaces over $X$, $x_0 \in X$ and $\xi_1, ...,\xi_n\in \IHH(x_0)$ be orthonormal. Then there is $U\in\mathfrak{U}(x_0)$ and $f_1, ... ,f_n\in\ICC(U)$ such that $f_1(x), ... , f_n(x)$ are orthonormal for all $x\in U$ and $f_j(x_0)=\xi_j$ for all $j=1, ... ,n$.
\end{proposition}
\begin{proof} (1) is a straightforward extension of \cite{dix}, \S 1, Remarque, p. 231: in fact, using the metric $d_{\IFF_x}$ from (\ref{f-d}), one finds that $f_1(x), ... , f_n(x)$ are linearly independent, if and only if
 $$
 \inf_{c\in\K^n, |c|=1}d_{\IFF_x}\Big(0,\sum_{k=1}^nc_kf_k(x)\Big)>0 ,
 $$
 and the latter positivity extends to a neighborhood of $x_0$ by continuity.\\
 (2) Part (1) combined with the Gram-Schmidt orthonormalization, see also \cite{god}, Proposition 8, p. 83.
\end{proof}

 \subsection*{Some general constructions}
 \begin{definition} Let $(\IVV,\ICC,p)$ be a continuous field of locally convex spaces over $X$, and let $\IWW_x\subset \IVV_x$ be a closed subspace for any $x\in X$. Then $(\IWW,p)$ is a field of locally convex spaces, and $\ICC\cap \Gamma(\IWW)$ the \emph{induced} continuous subspace. If it is saturated, we call 
  $(\IWW,\ICC\cap \Gamma(\IWW),p)$ a \emph{continuous subfield} of   $(\IVV,\ICC,p)$.
 \end{definition}
We note in passing, that in the above situation the closedness of $\IWW_x$ and the local uniform closedness of $\ICC$ imply that $\ICC\cap \Gamma(\IWW)$ is locally uniformly closed as well. The following extension of \cite{dix}, Proposition 11, p. 240, gives a satisfactory characterization of the required regularity of the subspaces in the Fr\'echet case:

\begin{proposition}\label{prop:subfield} Let $X$ be paracompact, $(\IEE,\ICC,p)$ be a continuous field of Fr\'echet spaces over $X$, for all $x\in X$ let $\IFF_x\subset \IEE_x$ be a closed subspace denote by $d_{\IEE_x}$ the canonical metric from \eqref{f-d} above. Then $(\IFF,\ICC\cap\Gamma(\IFF),p)$ is a continuous subfield, if and only if for all $f\in\ICC$ the function
$$
X\longrightarrow [0,\infty),\quad  x\longmapsto d_{\IEE_x}(f(x),\IFF_x)=\inf \set{d_{\IEE_x}(f(x),\zeta)}{\zeta\in\IFF_x}
$$
is upper semicontinuous. 
\end{proposition}

\begin{proof}
 This works like in \cite{dix}; for the calculations it is worthwhile to note that $d_{\IEE_x}(\lambda \xi, \lambda \zeta)\le d_{\IEE_x}(\xi,\zeta)$ for all $\lambda\in\K$ with $|\lambda|\le 1$.
\end{proof}

We next turn to the canonical way of passing from a seminormed space to a normed space by dividing out the kernel of the seminorm, in fact the procedure is more general:
\begin{proposition}\label{quotient}
 Let $(\IVV,\ICC, p)$ be a continuous field of locally convex spaces with index set $I$ over $X$, and for all $x\in X$ assume that
 $$
 \INN_x\subset\bigcap_{\iota\in I} \ker(p_\iota(x))
 $$ 
 is a closed subspace. Then, with an obvious notation, the triple $(\IVV/\INN,\ICC /\INN, \overline{p})$ becomes a continuous field of locally convex spaces over $X$, where $\overline{p}_\iota(x)([\xi]):=p_\iota(x)(\xi)$. 
\end{proposition}
This follows readily from the fact that the quotient mapping is onto and that 
$\overline{p}_\iota$ is well-defined. The following observation has on obvious extension to semimetric locally convex spaces, we will only use it for normed spaces and therefore state it in this less general setting to avoid additional notational complications.
\begin{proposition}\label{completion2} Let $(\IEE,\ICC)$ be a continuous
field of normed spaces. Denote by $\tilde{\IEE}_x$ the completion of $\IEE_x$ and consider
$\IEE_x\subset \tilde{\IEE}_x$, as well as $\ICC\subset\Gamma(\tilde{\IEE})$. We can then construct $\overline{\ICC}$ relative to $\tilde{\IEE}$, as in \ref{Lemma2.4}. Then $(\tilde{\IEE}, \overline{\ICC})$ is a continuous field of Banach spaces. If the $\tilde{\IEE}_x$ are pre-Hilbert spaces, $(\tilde{\IEE}, \overline{\ICC})$ is a continuous field of Hilbert spaces.
\end{proposition}
The assertions are obvious apart from the fact that $\overline{\ICC}$ is saturated. This, in turn was established in Lemma \ref{completion} above.
\subsection*{The net machine}

In the literature, notions are quite often formulated in terms of convergence of sequences rather than continuity. Of course, this can be rephrased in our framework, as maybe would not have been necessary to state:
\begin{remark}
Let $(X,\mathfrak{T})$ be a topological space.
 \begin{itemize}
  \item [(1)] A sequence $x=(x_n)_{n\in\N}$ converges to $
  x_\infty$ in $(X,\mathfrak{T})$ iff $\overline{x}:\N\cup\{\infty\}\to X$ is continuous, where $\N\cup\{\infty\}$ carries the natural topology, in that all points $n\in\N$ are isolated and a neighborhood base for $\infty$ is given by the sets $\set{n}{n\ge m}\cup\{\infty\}$, $m\in\N$.
  \item  [(2)] More generally, if $A$ is a directed set, we can define an analogous topology $\mathfrak{T}_A$ on $A\cup\{\infty\}$ by the neighborhood base $\set{\{\alpha\}}{\alpha\in A}\cup\set{\set{\beta}{\beta\ge\alpha}\cup\{\infty\}}{\alpha\in A}$ and the net $x=\snet{x}{\alpha}{A}$ converges to $x_\infty$ iff
  $\overline{x}:(A\cup\{\infty\},\mathfrak{T}_A)\to (X,\mathfrak{T})$ is continuous.
 \end{itemize}
\end{remark}
\begin{lemma}[Pullback of continuous fields]\label{pull}
Let $X$ and $Y$ be topological spaces and $\Phi:Y\to X$ be continuous. 
\begin{itemize}
    \item [(1)] If $(\IVV, p=(p_i;i\in I))$ is a field of locally convex spaces over $X$, then $\Phi^*\IVV_y:=\IVV_{\Phi(y)}$, $\Phi^*p_i(y):=p_i(\Phi(y))$ for $i\in I$, $y\in Y$ defines a field of locally convex spaces over $Y$ denoted by $(\Phi^*\IVV,\Phi^*p)$. 
    \item [(2)] If, furthermore, $\ICC$ is a continuous subspace for $(\IVV, p)$, then $\set{f\circ \Phi}{f\in \ICC}$ is a continuous subspace for $(\Phi^*\IVV,\Phi^*p)$.
    \item [(3)] If  $(\IVV, \ICC, p)$ is a continuous field, then so is $(\Phi^*\IVV,\Phi^*\ICC,\Phi^*p)$, where 
    $$\Phi^*\ICC:=\overline{\set{f\circ \Phi}{f\in \ICC}}
    $$ 
    is defined as in Lemma \ref{Lemma2.4}-(1).
\end{itemize}
\end{lemma}
We will frequently use a combination of the previous results for which we introduce appropriate shorthand  notation. To ease matters, nets in $X$ will always be denoted by $x=\snet{x}{\alpha}{A}$ as above, and we will assume (without loss of generality of course), that $\infty\notin A$. If  $x=\snet{x}{\alpha}{A}$ converges, its limit will be called $x_\infty$ and with $\overline{x}:A\cup\{\infty\}\to X$ we denote the continuous mapping introduced above. Then any [continuous] field $(\IVV,\ICC,p)$ over $X$ or subspace $\ICC$ defines, by pullback  a [continuous] field over $A\cup\{\infty\}$ which we denote as 
$(\IVV_{\overline{x}} ,\ICC_{\overline{x}} ,p_{\overline{x}})$ and call \emph{the field induced by} the net $x$. That this induced field inherits the properties from the original one will be used in arguments by referring to the net machine. We point out that, due to Proposition \ref{prop:properties of the topology} (3), a section $f\in\Gamma(\IVV_{\overline{x}})$ belongs to $\ICC_{\overline{x}}=\overline{x}^*\ICC$ precisely when $A\cup\{\infty\}\rightarrow \IVV_{\overline{x}}$ is continuous which, by the definition of the topology on $A\cup\{\infty\}$, is obviously equivalent to $\lim_\alpha f(\alpha)=f(\infty)$ in $\IVV$.

\section{Operators between continuous fields}\label{sec:opsbetweencontfields}
In what follows, we will consider [continuous] fields over the same topological space $X$. Let $(\IEE,\ICC_\IEE, p)$  and  $(\IFF,\ICC_\IFF, q)$ be two continuous fields of locally convex spaces over $X$ and define
$$
\elle{E}{F}:=\set{T}{T\in\ILL(\IEE_x,\IFF_x) \mbox{  for some  }\pi(T):=x},
$$
where $\ILL(\IEE_x,\IFF_x)$ denotes the set of continuous linear operators between the indicated spaces. Then there is a natural notion of continuity for sections $T$ of $\elle{E}{F}$, namely:
\begin{equation}\label{stronglycts}
T \mbox{ \emph{is strongly continuous} }:\Leftrightarrow Tf\in\ICC_\IFF\mbox{, whenever}f\in\ICC_\IEE ,
\end{equation}
where $Tf=T(\cdot)f(\cdot)$. This leads to a natural field of locally convex spaces of operators:
\begin{definition}\label{defls} Let $(\IEE,\ICC_\IEE, r)$ resp. $(\IFF,\ICC_\IFF, q)$ be continuous fields of locally convex spaces over $X$, with index set $A$ resp. $B$, and denote by
$$
\elless{E}{F}=\elle{E}{F}
$$
the family above, where the fibers  $\ILL(\IEE_x,\IFF_x)$ are equipped with the strong operator topology. For an arbitrary non-empty set $M$ we denote by $\operatorname{Fin}(M)$ the set of non-empty, finite subsets of $M$. Using $B$ and $\operatorname{Fin}\left(\ICC_\IEE\right)$ we let
$$
I:=\operatorname{Fin}\left(\ICC_\IEE\right)\times  B\ni (F,\beta)\mapsto p_{(F,\beta)}(\cdot)(T(\cdot)):= \max_{f\in F} q_\beta(\cdot)(T(\cdot)f(\cdot)), 
$$
defining a seminorm $p_i(x)$ on $\ILL(\IEE_x,\IFF_x)$ for every $i\in I$, $x\in X$. Moreover, we define
$$
\ICC^\s(\IEE,\IFF):=\set{T\mbox{ section of }\ILL^\s(\IEE,\IFF)}{T\mbox{ is strongly continuous} },
$$
and use the abbreviation $\ICC^\s$ when there is no risk of misinterpretation.
\end{definition}

\begin{lemma}\label{operator-field}
In the situation above we have
 \begin{itemize}
  \item[\textrm{(1)}] $(\elless{E}{F}, p)$ is a field of locally convex spaces with index set $I$ over $X$; its fibers carry the strong operator topology.
  \item[\textrm{(2)}] $\ICC^\s$ is a locally uniformly closed continuous subspace for  $(\elless{E}{F},p)$.
  \item[\textrm{(3)}] If $(\IGG,\ICC_\IGG,t)$ is another continuous field of locally convex spaces over $X$, for any $T\in\ICC^\s(\IEE,\IFF)$ and $S\in \ICC^\s(\IFF,\IGG)$ one has $ST\in \ICC^\s(\IEE,\IGG)$ for the fiberwise composition.
 \end{itemize}
\end{lemma}
\begin{proof}
 (1) For given $x\in X$,  the strong operator topology on $\ILL(\IEE_x,\IFF_x)$ is generated by the family of seminorms
 \[
 r_{\Xi,\beta}: \ILL(\IEE_x,\IFF_x) \ni S\mapsto \max_{\xi\in\Xi}q_\beta(S\xi)\in [0,\infty),\quad  \beta\in B, \Xi\in\operatorname{Fin}\left(\IEE_x\right),
 \]
 since the $q_\beta$ generate the topology of $\IFF_x$. As $\ICC_{\IEE}$ is saturated, for $\xi\in \IEE_x$, there is $ f_\xi\in \ICC_\IEE$ such that $f_\xi(x) = \xi$. Therefore, with $i= (\{f_\xi; \xi\in\Xi\},\beta)$ we get $p_i(x) = r_{\Xi,\beta}$  yielding the claim.\\
 (2) That $\ICC^s$ is continuous is obvious from the definitions. Let us check  local uniform closedness: 
 We start from a section $S$ in $\ILL(\IEE,\IFF)$ which satisfies condition (2.1) from Definition \ref{cflcs} for $\ICC^s$ and the family $p_i, i\in I$ of seminorms defined above.  Using the situation at hand this means that for any $x\in X$, $i\in I$ and $\epsilon>0$ we find $T\in\ICC^s, U\in \mathfrak{U}(x)$ such that
 \[
 p_i(y)\big( S(y ) - T(y) \big)\leq \epsilon\ \text{for all}\ y\in U.
 \]
We need to verify that $S\in \ICC^s$, that is $S f\in \ICC_\IFF$ for any $f\in\ICC_\IEE$ . \\
Let $\beta\in B$, set $i=(\{f\},\beta)$. As $T\in\ICC^s$, then  $h:=Tf\in \ICC_F$ and for $U$ as above we get
\begin{align*}
    p_i(y)\big( S(y ) - T(y) \big) &= q_\beta (y) \big( S(y)f(y) - T(y)f(y) \big)
    = q_\beta (y) \big( S(y)f(y) - h(y) \big) \\
    &\leq \epsilon \ \text{for all}\ y\in U.
\end{align*}
Thus  using the (LUC) property for $\ICC_\IFF$ with $h$ we get that $g:=Sf\in\ICC_\IFF$.\\
(3) Follows easily from the definition of $\ICC^\s(\IEE,\IGG)$.
\end{proof}

\begin{remark}\label{remstrong}
 \begin{itemize}
  \item [\textrm{(1)}] Note that in the above proof we only needed $\ICC_\IEE$ to be saturated and $\ICC_\IFF$ to be locally uniformly closed.
  \item[\textrm{(2)}] A special case covered by the considerations above is
  $$
  \IEE':= \elles(\IEE,X\times \K),
  $$
  where the fibers then carry the weak-*-topology. Here, the trivial field $X\times\K$ is equipped with the continuous subspace $C(X)$. We point out that in this special case one has 
  $$
  p=\set{p_F(x)}{F\in\operatorname{Fin}(\ICC_\IEE), x\in X}
  $$
  with
  $$p_F(\cdot)\left(u(\cdot)\right)=\max_{f\in F}|u(\cdot)f(\cdot)|,\quad u\in \IEE',$$
  and that
  $$\ICC^\s(\IEE,X\times \K)=\set{u\text{ section of }\IEE'}{\forall\,f\in\ICC_\IEE: u f=u(\cdot)f(\cdot)\in C(X)}.$$
  \item[\textrm{(3)}] Note that we did not claim that $(\elless{E}{F},\ICC^\s, (p_i; i\in I))$ is a continuous field of locally convex spaces in the general situation, since it is completely unclear whether $\ICC^\s$ is saturated.
  \end{itemize}
\end{remark}

\begin{definition}
    Let $(\IEE,\ICC_\IEE, p)$ be a continuous field of locally convex spaces over $X$ with index set $I$. 
    \begin{itemize}
        \item[\textrm{(1)}] We define
        $$
        \ICC^\w:=\ICC^\w_\IEE:=\set{f\in \Gamma(\IEE)}{\text{ for all $u\in\ICC^\s(\IEE, X\times\K)$ one has $u(\cdot)\left(f(\cdot)\right)\in C(X)$}}.
        $$
        Obviously, $\ICC^\w_\IEE$ is a subspace of $\Gamma(\IEE)$ which contains $\ICC_\IEE$ and is therefore saturated. Additionally, for $F\in\operatorname{Fin}(\ICC^s(\IEE,X\times\K))$ and $x\in X$
        $$q_F(x)(\zeta):=\max_{u\in F}|u(\pi(\zeta))(\zeta)|, \quad \zeta\in\IEE_x,$$
        defines a seminorm on $\IEE_x$ and for $f\in \ICC^\w_\IEE$ the function
        $$X\rightarrow\R, x\mapsto q_F(x)(f(x)):=\max_{u\in F}|u(x)\left(f(x)\right)|$$
        is continuous so that $\ICC^\w_\IEE$ is continuous for $\left(\IEE, q=(q_F; F\in\operatorname{Fin}(\ICC^s(\IEE,X\times \K)))\right)$. Therefore, by Lemma \ref{operator-field}, $(\IEE,\ICC^\w_\IEE, q)$ is a continuous field of locally convex spaces over $X$ which we abbreviate by $\IEE^\w$. Moreover, we then denote the topology on $\IEE$ induced by $\ICC^\w$ according to Proposition \ref{prop:properties of the topology} by $\mathfrak{T}^\w$ and call it the \emph{weak topology}. 
        
        We point out that under the hypothesis that $\set{u(x)}{u\in \ICC^\s(\IEE, X\times\K)}$ is weak-*-dense in $\IEE_x'$ for every $x\in X$ (which holds in particular whenever $\ICC^s(\IEE,X\times \K)$ is saturated for $\IEE'$), the family of seminorms $(q_F(x); F\in \operatorname{Fin}(\ICC^s(\IEE,X\times\K)))$ generates the weak topology on each fiber $\IEE_x$ of $\IEE$ which justifies the term weak topology. 

        In contrast to the weak topology $\mathfrak{T}^\w$, by the \emph{strong topology} we mean the topology $\mathfrak{T}^\s$ induced by $\ICC$ according to Proposition \ref{strong-top} above.
        \item[\textrm{(2)}] We say that $(\IEE,\ICC_\IEE, p)$ has \emph{locally equicontinuous duals}, if it satisfies the property
        \begin{itemize}
            \item[\textrm{(LED)}]
              \begin{gather}\label{led}
                \forall\,x_0\in X, u\in\IEE'_{x_0}\,\exists\,U\in\mathfrak{U}(x_0), i\in I, C>0\,\forall\,x\in U:\quad |u(x)(\zeta)|\leq C p_i(x)(\zeta)\quad (\zeta\in \IEE_x).
              \end{gather}
        \end{itemize}
    \end{itemize}
\end{definition}


\begin{proposition}\label{prop: properties of the weak topology}
    Let $(\IEE,\ICC_\IEE,p)$ and $(\IFF,\ICC_\IFF,r)$ be continuous fields of locally convex spaces over $X$.
    \begin{itemize}
        \item[(1)] If $(\IEE,\ICC_\IEE,p)$ is a continuous field of Hilbert spaces, then it has locally equicontinuous duals and $\ICC^\s(\IEE,X\times \K)$ is saturated in $\IEE'$. 
        \item[(2)] Let $T\in \elle{\IEE}{\IFF}$ be strongly continuous. Then $T f\in \ICC^\w_\IFF$ for every $f\in\ICC^\w_\IEE$.
        \item[(3)] Assume that $(\IEE,\ICC_\IEE,p)$ has locally equicontinuous duals. Then, $\mathfrak{T}^\w\subset\mathfrak{T}^\s$.
    \end{itemize}
\end{proposition}

\begin{proof}
    (1) follows from Riesz' representation theorem for Hilbert spaces combined with the fact that $\ICC_\IEE$ is saturated.

    In order to prove (2), we denote by $T'(x)$ the transposed of $T(x)$, $x\in X$. Then, for $v\in\ICC^s(\IFF, X\times \K)$ and $h\in \ICC_\IEE$ the function $X\rightarrow\K, x\mapsto T'(x)v(x)\left(h(x)\right)=v(x)\left(T(x)h(x)\right)$ is continuous because $Th\in\ICC_\IFF$. Hence, $T'v\in \ICC^s_\IEE$ which implies, for $f\in \ICC^\w_\IEE$, the continuity of the function $X\rightarrow \K, x\mapsto v(x)\left(T(x)f(x)\right)=T'(x)v(x)\left(f(x)\right)$ so that $Tf\in\ICC^\w_\IFF$.

    (3) Let $I$ be the index set of $p$. The following sets
    \[
    U^\w:=U^\w(V,g,F,\epsilon)=\set{\zeta\in\IEE}{\pi(\zeta)\in V, \max_{u\in F}|u(\pi(\zeta))\left(\zeta - g(\pi(\zeta))\right) |<\epsilon}
    \]
    form a base of $\mathfrak{T^\w}$, where $V\subset X$ is open $g\in \ICC^\w_\IEE$, $F\in\mbox{Fin}(\ICC^\s(\IEE,X\times\K))$ and $\epsilon >0$. Let $\xi\in U^\w$ and $x:=\pi(\xi)$; we need to find a neighborhood $U^s$ of $\xi$ with respect to the topology $\mathfrak{T}^\s$ which is contained in $U^\w$.
    
    Since $(\IEE,\ICC_\IEE,p_\IEE)$ has locally equicontinuous duals, for $u\in F$ there are $U_u\in \mathfrak{U}(x)$, $i_u\in I$, and $C_u>0$ so that $|u(y)\zeta|\leq C_u p_{i_u}(y)(\zeta)$ for each $y\in U_u$, $\zeta\in\IEE_y$, where without loss of generality we assume $U_u\subset V$. We pick $h\in\ICC_\IEE$ with $h(x)=g(x)$. 
    As $g\in\ICC^\w_\IEE$, $h\in\ICC_\IEE\subset\ICC^\w_\IEE$, and $u\in\ICC^\s(\IEE,X\times\K)$, the function
    \[
    y\mapsto \big|u(y)(g(y) - h(y) )\big|
    \]
    is continuous and vanishes in $x$, implying the existence of $W_u\subset U_u$, $W_u\in\mathfrak{U}(x)$ so that
    $$
    \big|u(y)\left( g(y) - h(y) \right)\big|<\frac{\epsilon}{2}\quad (y\in W_u).
    $$
    For $\zeta\in\pi^{-1}(W_u\cap U_u)$ it holds
    \begin{align*}
    \big|u(\pi(\zeta))\left(\zeta - g(\pi(\zeta)) \right)\big| & 
    \leq \big|u(\pi(\zeta))\left( \zeta - h(\pi(\zeta)) \right)\big|  + \big|u(\pi(\zeta))\left( g(\pi(\zeta)) - h(\pi(\zeta)) \right)\big|\nonumber \\
    &\leq C_u p_{i_u}(\pi(\zeta)\left(\zeta - h(\pi(\zeta))\right) + \frac{\epsilon}{2},
    \end{align*}
    so that
    $$U\left(W_u\cap U_u,h,i_u,\frac{\epsilon}{2 C_u}\right):=\set{\zeta\in\IEE}{\pi(\zeta)\in W_u\cap U_u, p_{i_u}(\pi(\zeta))\left(\zeta-h(\pi(\zeta))\right)<\frac{\epsilon}{2 C_u}}$$
    is a $\mathfrak{T}^\s$-neighborhood of $\xi$ such that
    $$\forall\,\zeta\in U\left(W_u\cap U_u,h,i_u,\frac{\epsilon}{2 C_u}\right):\,\big|u(\pi(\zeta))\left(\zeta - g(\pi(\zeta)) \right)\big|<\epsilon$$
    Hence,
    $$U^\s:=\bigcap_{u\in F}U\left(W_u\cap U_u,h,i_u,\frac{\epsilon}{2 C_u}\right)$$
    is a $\mathfrak{T}^\s$-neighborhood of $\xi$ with $U^\s\subset U^\w$, giving the claim.
\end{proof}

Even if we do not have a bona fide topology on the total space $\elle{E}{F}$, there is a natural ad hoc notion of 'strong convergence' for nets, used, e.g., in  \cite{kuwae}.  
\begin{definition}\label{def:strongoperator convergence}
 We say that a net $(T_\alpha)_{\alpha\in A}$ in $\elle{E}{F}$ \emph{converges meta strongly to} $T_\infty\in \elle{E}{F}$, if $\pi(T_\alpha)\to \pi(T_\infty)$, and, for every net $\snet{\zeta}{\alpha}{A\cup\{\infty\}}$ in $\IEE$ with $\pi(\zeta_\alpha)=\pi(T_\alpha)$ for all $\alpha\in A$:
 $$\lim_\alpha\zeta_\alpha =\zeta_\infty\mbox{ in }\IEE \Longrightarrow \lim_\alpha T_\alpha\zeta_\alpha = T_\infty\zeta_\infty \mbox{  in }\IFF.$$ 
 We then write 
 $$T_\alpha\xrightarrow{\ms} T_\infty$$
\end{definition}

Using the net $x=(\pi(T_\alpha))_{\alpha\in A}$, we see that, if $\pi(T_\alpha)\to \pi(T_\infty)$,
$$
T_\alpha\xrightarrow{\ms} T_\infty\quad\Longleftrightarrow\quad \snet{T}{\alpha}{A\cup\{\infty\}}\in \ICC^\s(\IEE_{\overline{x}},\IFF_{\overline{x}})
$$
by using the net machine.
We have the following elementary but useful observations:
\begin{lemma}\label{strongconvergence}
  Let $T\in \ICC^s(\IEE,\IFF)$ and $(T_\alpha)_{\alpha\in A\cup\{\infty\}}$ be a net in $\elle{E}{F}$.
  \begin{itemize}
      \item [(1)] $\| T\|$ is lower semicontinuous on $X$.
      \item[(2)] Assume that $(\elless{E}{F},\ICC^\s(\IEE,\IFF),p)$ is a continuous field of locally convex spaces. If $T_\alpha\xrightarrow{\ms} T_\infty$ then $\lim_\alpha T_\alpha=T_\infty$ with respect to the natural topology on $\elless{E}{F}$ induced by $\ICC^\s(\IEE,\IFF)$ according to Proposition \ref{prop:properties of the topology}.
\item [(3)] Consider the following conditions:
  \begin{itemize}
   \item[{\textrm (i)}]   $T_\alpha\xrightarrow{\ms} T_\infty$.
   \item[{\textrm (ii)}] There is $\alpha_0\in A$ such that $\sup_{\alpha\ge \alpha_0}\| T_\alpha\| <\infty$ and for every $\eta$ in a dense subset of $\IEE_{\pi(T_\infty)}$ there is a net $\snet{\eta}{\alpha}{A\cup\{\infty\}}$ in $\IEE$ with $\pi(\eta_\alpha)=\pi(T_\alpha)$ for all $\alpha\in A\cup\{\infty\}$ so that $\lim_\alpha\eta_\alpha =\eta_\infty$ and $\lim_\alpha T_\alpha\eta_\alpha = T_\infty\eta_\infty$.
  \end{itemize}
  Then $(ii) \Longrightarrow (i)$ and if $A$ admits a cofinal sequence $\snet{\alpha}{n}{\N}$, the converse holds.
  \item [(4)] $\| T_\infty\|\le \liminf\| T_\alpha\|$ whenever $T_\alpha\xrightarrow{\ms} T_\infty$.
   \end{itemize}
\end{lemma}
\begin{proof}
(1) As $T\in \ICC^\s(\IEE,\IFF)$, for $f\in\ICC_\IEE$ we get $Tf\in\ICC_\IFF$ and hence $\|T(\cdot)f(\cdot)\|\in C(X)$. Thus  $\| T\|$ is a supremum of continuous functions, hence lower semicontinuous.\\
(2) Let $B$ be the index set of the continuous field of locally convex spaces $(\IFF,\ICC_\IFF, q)$. Since for $f\in\ICC_\IEE$ and $\Theta\in \ICC^\s(\IEE,\IFF)$ we have $\Theta f \in\ICC_\IFF$, it follows from Remark \ref{rem:natural topology is an initial topology} that
$$
\IFF\longrightarrow\R,\quad  \eta\longmapsto \max_{f\in F}q_\beta(\pi(\eta))\left(\eta-\Theta(\pi(\eta))f(\pi(\eta)\right)
$$
is continuous for every $\beta\in B$ and $F\subset\operatorname{Fin}(\ICC_\IEE)$. Therefore, as by hypothesis $\lim_\alpha T_\alpha f(\pi(T_\alpha))=T_\infty f(\pi(T_\infty))$ for any $f\in\ICC_\IEE$, we conclude
\begin{gather*}
    \lim_\alpha \max_{f\in F}q_\beta\left(\pi(T_\alpha)\right)\left(T_\alpha f(\pi(T_\alpha))-\Theta(\pi(T_\alpha))f(\pi(T_\alpha)) \right)\\=\max_{f\in F}q_\beta\left(\pi(T_\infty)\right)\left(T_\infty f(\pi(T_\infty))-\Theta(\pi(T_\infty))f(\pi(T_\infty)) \right)
\end{gather*}
for arbitrary $\beta\in B$, $F\in\operatorname{Fin}(\ICC_\IEE)$, and $\Theta\in\ICC^\s(\IEE,\IFF)$. Therefore, for every $Y\in\mathfrak{U}(\pi(T_\infty))$, $\Theta\in \ICC^\s(\IEE,\IFF)$, $(F,\beta)\in\operatorname{Fin}(\ICC_\IEE)\times B$, and $\epsilon>0$, taking into account the continuity of $\pi$ as well as
$$T_\infty\in U(Y,\Theta, (F,\beta), \epsilon):=\set{S\in \elless{\IEE}{\IFF}}{ \pi(S)\in Y, \max_{f\in F}q_\beta(\pi(S))\left(S-\Theta(\pi(S))f(\pi(S))\right)}$$
we obtain the existence of $\alpha_0$ such that $T_\alpha\in U(Y,\Theta, (F,\beta), \epsilon)$ whenever $\alpha\geq \alpha_0$. By assumption, $(\elless{E}{F},\ICC^\s(\IEE,\IFF),p)$ is a continuous field of locally convex spaces so that the sets $U(Y,\Theta, (F,\beta), \epsilon)$ from above form a neighborhood basis of $T_\infty$ with respect to the natural topology on $\elless{E}{F}$ which proves (2). 

(3) $(ii)\Rightarrow (i)$:

Let $\xi\in\IEE_{\pi(T_\infty)}$ and let $\xi_\alpha\in \IEE_{\pi(T_\alpha)}$ be such that $\lim_\alpha\xi_\alpha = \xi_\infty$. We have to show that $\lim_\alpha T_\alpha\xi_\alpha = T_\infty\xi_\infty$. For this, let $U(Y,g,\delta)$ be a neighbourhood of $T_\infty\xi_\infty$, $Y\subset X$ open, $g\in \ICC_\IFF$, $\delta>0$.

With $\alpha_0$ from the hypotheses we set $C:=\sup_{\alpha\geq \alpha_0}\|T_\alpha\|$, and we denote the dense subspace of the $\eta$ from the hypotheses by $E$. Let $U(Y,f,\epsilon)$ be another neighbourhood of $T_\infty\xi_\infty$ such that $U(Y,f,3 \epsilon)\subset U(Y,g,\delta)$. We find $\eta\in E$ such that 
$$
\|\eta - \xi_\infty\|\leq \epsilon/(1+C+\|T_\infty\|),
$$ 
and then
\begin{equation}\label{eq:characterization of strong operator convergence 1}
   \|T_\infty\eta - f(\pi(T_\infty\eta))\| = \|T_\infty\eta_\infty - f(\pi(\xi_\infty))\| \leq \| T_\infty\eta - T_\infty\xi_\infty\| + \|T_\infty\xi_\infty - f(\pi(\xi_\infty))\| \leq 2\epsilon. 
\end{equation}
As $\eta\in E$, by hypothesis there is a net $(\eta_\alpha)$ such that $\pi(\eta_\alpha) = \pi(T_\alpha)$, $\lim_\alpha\eta_\alpha =\eta$ and $\lim_\alpha T_\alpha\eta_\alpha = T_\infty\eta$. In order to conclude $T_\alpha\xi_\alpha\in U(Y,f,\epsilon)$ for every $\alpha\geq \alpha_1$ with suitable $\alpha_1$, we observe that on the one hand we have 
$$
\pi(T_\alpha\xi_\alpha) = \pi(\xi_\alpha) = \pi(\eta_\alpha)\to \pi(\eta)= \pi(T_\infty\eta) = \pi(T_\infty\xi_\infty),
$$
by continuity of $\pi$ and then $T_\alpha\xi_\alpha\in Y, \alpha\geq\tilde{\alpha}_1$ for some $\tilde{\alpha}_1$. On the other hand, for $\alpha\geq \alpha_0$ we have
\begin{align*}
\| T_\alpha\xi_\alpha - f(\pi(T_\alpha\xi_\alpha)) \| & \leq \| T_\alpha\xi_\alpha - T_\alpha\eta_\alpha \|
 + \| T_\alpha\eta_\alpha - f(\pi(T_\alpha\eta_\alpha)) \|\\
 & \leq C \|\xi_\alpha - \eta_\alpha\| + \| T_\alpha\eta_\alpha - f(\pi(T_\alpha\eta_\alpha)) \|.
\end{align*}
As $\lim_\alpha \xi_\alpha=\xi_\infty$, $\lim_\alpha \eta_\alpha=\eta$, and $\lim_\alpha T_\alpha\eta_\alpha=T_\infty \eta$, passing to the limit in the previous inequality in combination with the continuity of $f\circ \pi$, Proposition \ref{prop:properties of the topology} (2), and \eqref{eq:characterization of strong operator convergence 1} now gives
\[
\lim\sup_\alpha\| T_\alpha\xi_\alpha - f(\pi(T_\alpha\xi_\alpha)) \| \leq C \|\xi_\infty - \eta\| 
+ \| T_\infty\eta - f(\pi(T_\infty\eta)) \|<3\epsilon.
\]
Hence, there is $\alpha_1$ such that $T_\alpha\xi_\alpha\in U(Y,f,3\epsilon)\subset U(Y,g,\delta)$ whenever $\alpha\geq \alpha_1$ implying $\lim_\alpha T_\alpha\xi_\alpha = T_\infty\xi_\infty$, as desired.\\

$(i)\Rightarrow (ii)$. Because $\ICC_\IEE$ is saturated and because for every $f\in \ICC_\IEE$ the net $\eta_\alpha:=f(\pi(T_\alpha))$ satisfies $\lim_\alpha \eta_\alpha= f(\pi(T_\infty))$, we only have to show $\sup_{\alpha\ge \alpha_0}\| T_\alpha\| <\infty$. Assume the contrary. Then, for any $n\in\N$ there is $\beta(n)\in A$ so that $\beta(n)\ge\alpha_n$ and $\xi_{\beta(n)}\in \IEE_{\pi(\xi_{\beta(n)})}$ with $\|\xi_{\beta_n}\| \le \frac{1}{n}$ and $ \|T_{\beta(n)} \xi_{\beta(n)}\| = 1$. Let $\eta_\alpha:= \xi_{\beta(n)}$, for $\alpha=\beta(n)$ and $0$ else. By construction, $\eta_{\alpha}\to 0$ in $\IEE$, whereas $\|T_{\beta(n)} \eta_{\beta(n)}\| = 1$, which contradicts (i).

(4) is easy to prove with the help of Proposition \ref{prop:properties of the topology} (2).
\end{proof}

In \cite{stummel1}, (6), p. 53 the analogous equivalence in (2) is shown, Stummel speaks of 'Stabilität' for the uniform boundedness (as he is dealing with sequences) and 'Konsistenz' for the second condition in (ii) above.

Maybe, a word on notation is in order: we introduced $\xrightarrow{\ms}$ for operators since, following our general convention, $\xrightarrow{\s}$ is reserved for convergence w.r.t. the natural topology coming from $\ICC^\s$, once the latter induces a continuous field. By Example \ref{ex:weaknotmeta} below, $\ms$-convergence is in general strictly stronger than $\s$-convergence. See also Corollary \ref{prop:notopology} that shows that on a single infinite dimensional Hilbert space, $\ms$-convergence isn't even induced by a topology.

As nice as nets may be, they are technically more complicated to handle. Moreover, when working with spaces of integrable functions, sequences allow the application of limit theorems from integration theory. From this perspective, the following observation is quite useful.

\begin{corollary}
    \label{cor:seqcontinuity}
    Let $X$ be first countable and $(\IEE,\ICC_\IEE)$, $(\IFF,\ICC_\IFF)$ be continuous fields of Fr\'echet spaces. For $T\in\Gamma(\elle{\IEE}{\IFF})$ the following are equivalent:
    \begin{itemize}
        \item [(i)] $T\in\ICC^\s(\IEE,\IFF)$.
        \item [(ii)] For every $f\in\ICC_\IEE$, $x_\infty\in X$ and every sequence $\snet{x}{n}{\N}$ in $X$ so that $x_n\to x_\infty$:
        $$
        T_{x_n}f(x_n)\xrightarrow{\s}T_{x_\infty}f(x_\infty).
        $$
    \end{itemize}
    If, furthermore, $(\IEE,\ICC_\IEE)$, $(\IFF,\ICC_\IFF)$ are continuous fields of Banach spaces and $\| T\|$ is locally bounded, then \emph{(i)} and \emph{(ii)} are in turn equivalent to
    \begin{itemize}
        \item [(iii)] For every $x_\infty\in X$ and every sequence $\snet{x}{n}{\N}$ in $X$ so that $x_n\to x_\infty$:
        $$
        T_{x_n} \xrightarrow{\ms}T_{x_\infty}.
        $$
    \end{itemize}    
\end{corollary}
\begin{proof}
    For the first part, simply note that the continuity of $Tf$ is equivalent to sequential continuity by virtue of Corollary \ref{cor:countable}.

    For the second part, it suffices to observe that for any $\snet{\zeta}{n}{\N}$ in $\IEE$ with $\pi(\zeta_n)=x_n$ and $\zeta_n\xrightarrow{\s}\zeta_\infty$ and any $f\in\ICC_\IEE$ with $f(x_\infty)=\zeta_\infty$ it follows that $\| f(x_n)-\zeta_n\|\to 0$ as $n\to\infty$.
\end{proof}

We will now show that $\ICC^\s$ is saturated in the Hilbert case under rather mild conditions on the underlying topological space.
\begin{lemma}\label{findim}
Let $(\IHH,\ICC_\IHH)$ and $(\IKK,\ICC_\IKK)$ be continuous fields of Hilbert spaces.
 Assume that for $x_0\in X$ there are orthonormal systems $(\xi_1, ..., \xi_n)$ in $\IHH_{x_0}$ and $(\eta_1, ..., \eta_n)$ in $\IKK_{x_0}$. Then there is an open $U\in \mathfrak{U}(x_0)$ and $P\in\ICC^\s(\IHH|_U,\IKK|_U)$ such that $P(x_0)\xi_k=\eta_k$ and $P(x)$ is a partial isometry from an $n$-dimensional subspace of $\IHH_x$ onto an an $n$-dimensional subspace of $\IKK_x$ for every $x\in U$.
\end{lemma}
\begin{proof}
 This follows from the existence of continuous orthonormal frames established in Proposition \ref{ON-Frames}-(2). In fact, we find an open set $U\in\mathfrak{U}(x_0)$ and $f_j\in\ICC(\IHH|_U), g_j\in\ICC(\IKK|_U), j=1,\cdots,n$ such that  $f_j(x_0) = \xi_j, g_j(x_0) = \eta_j$, so that 
 \[
\langle f_j(x), f_k(x) \rangle = \langle g_j(x), g_k(x) \rangle = \delta_{jk}\ \text{for all}\  x\in U.
 \]
 We now take 
 \[
P(x) = \sum_{j=1}^n \langle \cdot, f_j(x) \rangle g_j(x), x\in U.
 \]
 It follows that $P(x)$ is a unitary operator from ${\rm {lin}}\{ f_1(x),\cdots,f_n(x)\}$ into ${\rm {lin}}\{ g_1(x),\cdots,g_n(x)\}$, $x\in U$.\\
 To prove strong continuity of  $P$, pick $f\in \ICC(\IHH|_U)$ , then
 \[
(Pf)(x) = \sum_{j=1}^n \langle f(x), f_j(x) \rangle g_j(x), x\in U.
 \]
 By choice, $f_j \in \ICC(\IHH|_U), g_j\in\ICC(\IKK|_U)$, moreover by polarization the scalar functions 
 $$
 \langle f(\cdot), f_j(\cdot) \rangle \quad\text{are in $C(U)$},
 $$
 which leads to $Pf\in\ICC^s( \IHH|_U, \IKK|_U )$ since the latter is a module with respect to $C(U)$.
\end{proof}
\begin{theorem}\label{unitary}
Let $X$ be a locally compact Hausdorff space or a normal space, $x_0\in X$, $(\IHH,\ICC_\IHH)$ and $(\IKK,\ICC_\IKK)$ be continuous fields of Hilbert spaces over $X$. Assume that $\net{\xi}{n}{\N}$ and $\net{\eta}{n}{\N}$ are orthonormal bases for $\IHH_{x_0}$ and $\IKK_{x_0}$, respectively. Then there is $P\in\ICC^\s(\IHH,\IKK)$ so that $P(x_0)\xi_k=\eta_k$ for $k\in\N$. 
\end{theorem}
Note that, in particular, $P(x_0)$ is unitary.
\begin{proof}
 Using Lemma \ref{findim} for fixed $n$, we find $U_n\in\mathfrak{U}(x_0)$ and $\tilde{f}_j\in\ICC_\IHH(U_n), \tilde{g}_j\in\ICC_\IKK(U_n)$, so that $\tilde{f_j}, \tilde{g_j}, j=1,\cdots,n$ are orthonormal frames on $U_n$. Making those sets smaller, if necessary we can  arrange that $\overline{U}_1\subset U_0, \overline{U}_1$ compact and  $\overline{U}_{n+1} \subset U_n$. Let us choose further cut-offs $\varphi_n\in C_c(X), 1_{U_n}\leq \varphi_n \leq 1_{U_{n-1}}$ so that $f_j:=\tilde{f}_j\varphi_j\in \ICC_\IHH, g_j:=\tilde{g}_j\varphi_j\in \ICC_\IKK$, giving
 \[
 1_{U_n} \leq \|f_j\| , \|g_j\|\leq 1_{U_{n-1}}\  \text{for}\ 1\leq j\leq n, \langle f_j,f_k \rangle= \langle g_j,g_k \rangle =\delta_{jk}\ \text{for}\ k\neq j.
 \]
 Moreover we can assume that $\cap_{n} U_n =\{x_0\}$.\\
 We define
 \[
 P(x) := \sum_{j=1}^\infty \langle \cdot, f_j(x) \rangle g_j(x),
 \]
 the convergence being in the strong operator topology. Clearly $P(x_0)$ is an isometry.\\
 For any $x\in X\setminus\{x_0\}$, there is $N\in\N$ such that $x\in X\setminus U_N$ and hence 
 \[
P(x) = \sum_{j=1}^N \langle \cdot, f_j(x) \rangle g_j(x),
 \]
 Arguing as in the proof of Lemma \ref{findim} we conclude that $P(\cdot)$ is strongly continuous on the open set $X\setminus\{x_0\}$.\\
 It remains to prove continuity in $x_0$. To this end we will use Lemma \ref{strongconvergence}. Let $(x_\alpha)_{\alpha\in A}$ be a net in $X\setminus\{x_0\}$ converging to $x_0$ and observe that the subspace $E_0:=\big\{{\rm{lin}}\{f_1(x_0),\cdots,f_n(x_0)\}, n\in\N\big\}$ is dense in $\IHH_{x_0}$. Fix $\xi\in E_0, f\in\ICC_\IHH$ such that $f(x_0) = \xi$ and write 
 \[
 \xi =\sum_{j=1}^n \langle f(x_0),f_j(x_0) \rangle f_j(x_0), \xi_\alpha :=\sum_{j=1}^n \langle f(x_\alpha),f_j(x_\alpha) \rangle f_j(x_\alpha) \in\IHH_{x_\alpha}
 \]
 Then $\pi(\xi_\alpha) = \pi(P(x_\alpha))$ for all $\alpha$ and $\xi_\alpha\xrightarrow{\s} \xi$, since $f,f_j$ are in $\ICC_\IHH$ for all $j$.\\
 Let $\eta\in \IHH_{x_\alpha}$, since the families $f_j, g_j$ are orthonormal we obtain
 \[
\|P(x_\alpha)\eta\|^2 = \sum_{j=1}^\infty \big( \langle \eta,f_j(x_\alpha) \rangle \big)^2 \leq \|\eta\|^2 ,
 \]
 leading to $\sup_{\alpha\in A}  \|P(x_\alpha)\|^2\leq 1 $. Moreover, 
 \[
P(x_\alpha)\xi_\alpha =\sum_{j=1}^n \langle f(x_\alpha),f_j(x_\alpha) \rangle g_j(x_\alpha) \xrightarrow{\s} \sum_{j=1}^n \langle f(x_0),f_j(x_0) \rangle f_j(x_0)
= P(x_0)\xi.
 \]
Owing to Lemma \ref{strongconvergence}, the above considerations lead to strong convergence of $P(x_\alpha)$ towards $P(x_0)$. Since $P(\cdot)$ is strongly continuous on $X\setminus\{x_0\}$, using (LUC') property we conclude that $P$ is continuous at $x_0$ and hence strongly continuous on $X$. 
 \end{proof}
 It is clear from the proof that we can choose $P(x_0) = \operatorname{Id}_{\IHH_{x_0}}$, the identity on $\IHH_{x_0}$, provided $\IHH_{x_0} = \IKK_{x_0}$.\\
This was the main step for the proof of the following result: 
\begin{theorem}\label{niceone}
Let $X$ be a locally compact Hausdorff space, or a normal space, and $(\IHH,\ICC_\IHH)$ and $(\IKK,\ICC_\IKK)$ be continuous fields of Hilbert spaces over $X$. Assume that $\IHH_x$ and $\IKK_x $ are separable for every $x\in X$. Then $\ICC^\s(\IHH,\IKK)$ is saturated. 
\end{theorem}
\begin{proof}
Let  $x_0\in X$, and $T_0\in \ILL(\IHH_{x_0},\IKK_{x_0})$. We have to find  $T\in \ICC^s(\IHH,\IKK)$ such that $T(x_0) = T_0$. By the previous theorem we find $P\in\ICC^s(\IHH, X\times \IHH_{x_0})$ such that $P(x_0) = \operatorname{Id}_{\IHH_{x_0}}$.\\
Again by the previous theorem we find $Q\in\ICC^s(X\times \IHH_{x_0}, \IKK)$  such that $Q(x_0) = \operatorname{Id}_{\IKK_{x_0}}$. Obviously the constant field $T_0(\cdot):=T_0: X\times \IHH_{x_0}\to X\times \IKK_{x_0}$ is continuous. Thus, by Lemma \ref{operator-field},   $Q\circ T_0(\cdot)\circ P$ does the job.
\end{proof}
  Using Theorem \ref{niceone} in conjunction with  Lemma \ref{operator-field}, we obtain:
\begin{corollary}\label{Ls}
 Let $X$ be a locally compact Hausdorff space, or a normal space, and $(\IHH,\ICC_\IHH)$ and $(\IKK,\ICC_\IKK)$ be continuous fields of Hilbert spaces over $X$. Assume that $\IHH_x$ and $\IKK_x$ are separable for every $x\in X$. Then $(\elless{H}{K},\ICC^\s(\IHH,\IKK),p)$ is a continuous field of locally convex spaces, where $p$ is given in {\emph{Definition \ref{defls}}} above. 
\end{corollary}

\section{Weak and strong convergence on the total space of a continuous field of Hilbert spaces}
\label{sec:weakandstrong}

In this section, we fix a continuous field of Hilbert spaces $(\IHH,\ICC)$ over $X$. We will see that on the total space $\IHH$ the interplay between the weak and strong topology works pretty much like in Hilbert spaces. To fix notation, we state: 

\begin{definition} Let $(\zeta_\alpha)_{\alpha\in A}$ be a net in $\IHH$ and let $\zeta_\infty\in \IHH$.
\begin{itemize}
    \item [(1)] We say that $\snet{\zeta}{\alpha}{A}$ \emph{converges strongly to $\zeta_\infty$}, if $\lim_{\alpha}\zeta_\alpha=\zeta_\infty$ in $(\IHH,\mathfrak{T}^\s)$. We then write $$\zeta_\alpha\xrightarrow{\s} \zeta_\infty$$

   \item [(2)] We say that $\snet{\zeta}{\alpha}{A}$ \emph{converges weakly to $\zeta_\infty$}, if $\lim_{\alpha}\zeta_\alpha=\zeta_\infty$ in $(\IHH,\mathfrak{T}^\w)$. We then write 
$$\zeta_\alpha\xrightarrow{\w} \zeta_\infty$$ 

   \item [(3)] We say that $(\zeta_\alpha)$ \emph{converges meta weakly to $\zeta_\infty$}, if for every net $\snet{\eta}{\alpha}{A}$ in $\IHH$ with $\pi(\eta_\alpha)=\pi(\zeta_\alpha)$ for all $\alpha\in A$,
$$\eta_\alpha\xrightarrow{\s} \eta_\infty\ \quad\Longrightarrow\quad \left\langle\zeta_\alpha,\eta_\alpha\right\rangle\to  \left\langle\zeta_\infty,\eta_\infty\right\rangle .$$
We then write 
$$\zeta_\alpha\xrightarrow{\m\w} \zeta_\infty$$
\end{itemize}
\end{definition}
The convergence we call meta weakly is called weak convergence, e.g. in the fundamental paper of Kuwae and Shioya, \cite{kuwae}. However, in the case of a single infinite dimensional space it does not coincide with weak convergence, the crucial difference being that a weakly convergent net need not be bounded. Of course, any weakly convergent sequence is bounded and therefore meta weakly convergent, as is readily seen. We record some observations that are easy to prove:
\begin{proposition}
    For an eventually bounded net  $(\zeta_\alpha)_{\alpha\in A}$ in $\IHH$:
    $$
  \zeta_\alpha\xrightarrow{\w} \zeta_\infty\Longrightarrow \zeta_\alpha\xrightarrow{\m\w} \zeta_\infty$$
\end{proposition} 

\begin{proof}
The claim is an immediate consequence of Lemma \ref{strongconvergence}.
\end{proof}

In conjunction with Lemma \ref{strongconvergence}-(4), we immediately get the lower semicontinuity of the norm w.r.t. weak convergence:

\begin{corollary}\label{cor:lscnorm} 
$\|\cdot\|:\IHH^\w\to\R$ is lower semicontinuous; in particular, $\| f\|:X\to\R$ is lower semicontinuous for every $f\in\ICC^\w$.
\end{corollary}

As indicated above, weak convergence does not imply meta weak convergence. We include an example for the readers convenience:

\begin{example}\label{ex:weaknotmeta}
Let $\IHH_0$ be an infinite dimensional Hilbert space. Denote by $\mathfrak{U}^w(0)$ the filter of neighborhoods of $0$ in the weak topology; note that for any $U\in \mathfrak{U}^w(0)$ there is an element $x_U\in U$ with $\| x_U\| =1$ and such that $\R x_U\subset U$ (just consider the standard neighborhood basis and use the fact that $\IHH_0$ is infinite dimensional).
\begin{itemize}
    \item [(1)] Define $A:= \mathfrak{U}^w(0) \times \N$ with the partial order
    $$
    (U,k)\le (V,n) :\Longleftrightarrow V\subset U \mbox{  and  } n\ge k 
    $$
    and a net $(x_\alpha)_{\alpha\in A}$ by $x_{(U,k)}:= k\cdot x_U$. Then
    $$
    x_\alpha\xrightarrow{\w} 0 .
    $$
    \item [(2)] For   $y_{(U,k)}:= \frac{1}{k}\cdot x_U$, we get 
    $$
    y_\alpha\xrightarrow{\s} 0 \mbox{  and  }\langle x_\alpha,y_\alpha\rangle =1
    $$
    so that, consequently, $(x_\alpha)_{\alpha\in A}$ does not converge meta weakly.
\end{itemize}
\end{example}
In view of the following result it is clear that meta weak convergence is strictly stronger than weak convergence. However, we included the preceding example since it plays a crucial role in the proof of that result.

We should note at this point, that the proof of Lemma 2.3 in \cite{kuwae} is not completely convincing, the reason being an erroneous interpretation of subnets. Even if this is unnecessary we recall the concept of a subnet (due to Kelley and Willard \cite{kelley}): a \emph{subnet} of a net $\snet{\zeta}{\alpha}{A}$ is given by a mapping  $B\xrightarrow{\varphi} A$ that is order preserving and has cofinal image, i.e., for every $\alpha\in A$ there is $\beta\in B$ with $\varphi(\beta)\ge \alpha$.

We should also point out, that meta weak convergence does not induce a topology in general, even not on a single infinite dimensional Hilbert space. Therefore, Definition 2.5 in \cite{kuwae} seems a little bit too optimistic:

\begin{proposition}\label{prop:notopology}
Let $\IHH_0$ be an infinite dimensional Hilbert space. 
\begin{itemize}
    \item [(1)]  A net $\snet{\xi}{\alpha}{A}$ with $\|\xi_\alpha\|\to \infty$ is not meta weakly convergent. 
    \item [(2)] 
There is no topology on $\IHH_0$ such that the convergent nets with respect to that topology are precisely the meta weakly convergent nets.  
\end{itemize}
\end{proposition}
\begin{proof}
(1) is proved in the exact same way as (2) in the preceding example.

(2) Assume, for a contradiction, that there is a topology that induces meta weak convergence and denote, by $\mathfrak{V}$ the corresponding filter of neighborhoods of $0$. First observe that the corresponding topology is linear, due to the evident 'linearity' of meta weak convergence.

Case (A): There is a norm-bounded $V\in\mathfrak{V}$. Then, since all $rV$, $r>0$ belong to $\mathfrak{V}$ as well, the identity $(\IHH_0,\mathfrak{V})\rightarrow \IHH_0$ is continuous, where, as range space, $\IHH_0$ is equipped with the norm topology. On the other hand, for a sequence $(x_n)_{n\in\N}$ in $\IHH_0$ converging in the norm topology to $x_\infty$, due to the norm boundedness of $\{x_n-x_\infty| n\in\N\}$, $(x_n)_{n\in\N}$ converges meta weakly to $x_\infty$. Hence, $id:\IHH_0\rightarrow (\IHH_0,\mathfrak{V})$ is sequentially continuous, thus, by \cite[Proposition 1.6.15]{Engelking} continuous, resulting in the equality of the norm topology and $\mathfrak{V}$. However, as $\IHH_0$ is infinite dimensional, there is an orthonormal sequence $(e_n)_{n\in\N}$ which is readily seen to converge meta weakly to zero which gives the desired contradiction.

Case (B): Every $V\in\mathfrak{V}$ is unbounded. Then a construction as in (1) of the preceding Example gives a net $\snet{\xi}{\alpha}{A}$ with $\|\xi_\alpha\|\to \infty$ that converges meta weakly to $0$.  
\end{proof}
We state an immediate consequence of the arguments in the preceding proof:

\begin{corollary}
    Let $(\IEE_0,\|\cdot\|)$ be a normed space and let $\mathfrak{T}$ be a linear topology $\IEE_0$ that is weaker than the norm topology $\mathfrak{T}_{\|\cdot\|}$. Then either $\mathfrak{T}=\mathfrak{T}_{\|\cdot\|}$ or there exists a net $\snet{x}{\alpha}{A}$ in $\IEE_0$  with $\| x_\alpha\|\to\infty$ and $x_\alpha\xrightarrow{\mathfrak{T}}0$.
\end{corollary}

\begin{proposition}\label{prop:bddgivesmw} Let $(\zeta_\alpha)_{\alpha\in A}$ be eventually bounded and $\set{\pi(\zeta_\alpha)}{\alpha\in A}$ be relatively compact in $X$. Then there exists a
meta weakly convergent subnet of $(\zeta_\alpha)_{\alpha\in A}$.    
\end{proposition}

\begin{proof} Passing to a subnet if necessary, we may assume that the $\pi(\zeta_\alpha)$ converge to some $x\in X$ and that, moreover $\| \zeta_\alpha\|\le 1$ for all $\alpha$ (by scaling). To use Tychonoff's theorem we put
$$
\INN :=\set{\snet{\eta}{\alpha}{A}}{\pi(\eta_\alpha)=\pi(\zeta_\alpha),\, \eta_\alpha \xrightarrow{\s}\eta_\infty\mbox{ and }\|\eta_\alpha\|\le \|\eta_\infty\|\,(\alpha\in A)}.
$$
By construction for any $\overline{\eta}\in\INN$, and $\alpha\in A$
$$
\overline{\zeta}_\alpha:=\left(\langle \zeta_\alpha,\eta_\alpha\rangle\right)_{\overline{\eta}\in\INN}
\in\Pi_{\overline{\eta}\in\INN}\overline{B_{\|\eta\|}(0)},
$$
where $\overline{B_r(0)}\subset\K$ denotes the closed ball of radius $r$ centered at $0$. Tychonoff's theorem gives a subnet  $\snet{\overline{\zeta}}{\beta}{B}$ that converges in the product topology, which means that, for any $\overline{\eta}\in\INN$, $\lim_{\beta\in B}\langle \zeta_\beta,\eta_\beta\rangle$ exists. Since for any $\eta_\infty\in \IHH_x$ there is $f\in\ICC$ with $f(x)=\eta_\infty$, setting $$\eta_\alpha:=\|f(x)\|_x\left(\max\{1,\|f(\pi(\zeta_\alpha))\|_{\pi(\zeta_\alpha)}\}\right)^{-1}f(\pi(\zeta_\alpha))\quad(\alpha\in A),$$
we obtain $(\eta_\alpha)_{\alpha\in A}\in\INN$ with $\eta_\alpha \xrightarrow{\s}\eta_\infty$.

It is easily seen that for $\eta_\infty\in \IHH_x$ the limit $\lim_{\beta\in B}\langle \zeta_\beta,\eta_\beta\rangle$ is independent of the choice of $(\eta_\alpha)_{\alpha\in A}$ in $\INN$ with $\eta_\alpha \xrightarrow{\s}\eta_\infty$ and that the functional
$\eta_\infty\mapsto \lim_{\beta\in B}\langle \zeta_\beta,\eta_\beta\rangle$ is linear and continuous. Consequently, there is $\zeta\in\IHH_x$ so that $\lim_{\beta\in B}\langle \zeta_\beta,\eta_\beta\rangle=\langle \zeta,\eta\rangle$. An appeal to Lemma \ref{strongconvergence} gives that   $\lim_{\beta\in B}\langle \zeta_\beta,\xi_\beta\rangle=\langle \zeta,\xi\rangle$ for all strongly convergent nets $\snet{\xi}{\beta}{B}$ as was to be proven.
\end{proof}

\begin{proposition}\label{prop:mwandnormgivesstrong}
\begin{itemize}
\item[(1)]
If  $(\zeta_{\alpha})_{\alpha\in A}$ is a net in $\IHH$ with $\zeta_\alpha\xrightarrow{\w} \zeta_\infty$ and $\left\|\zeta_\alpha\right\|\to  \left\|\zeta_\infty\right\|$, then $\zeta_\alpha\xrightarrow{\s} \zeta_\infty$.
\item[(2)] If $\zeta_\infty\in \IHH$ and  $(\zeta_{\alpha})_{\alpha\in A}$ is a net in $\IHH$ with the property that for every meta weakly convergent net $\snet{\xi}{\alpha}{A}$ with $\pi(\zeta_\alpha)=\pi(\xi_\alpha)$, $\langle \zeta_\alpha,\xi_\alpha\rangle\to \langle \zeta_\infty,\xi_\infty\rangle$, then $\zeta_\alpha\xrightarrow{\s} \zeta_\infty$.  
\end{itemize}
    \end{proposition}

\begin{proof} The proof for the Hilbert space case extends verbatim to our setting here.
\end{proof}
\iffalse
\texttt{... gerne hätte ich noch...}
\begin{theorem}
 Let $Y\subset X$ and $C_y\subset\IHH_y$ be convex for all $y\in Y$. Then
 $$
 \overline{\bigcup_{y\in Y}C_y}^\w=\overline{\bigcup_{y\in Y}C_y}^\s .
 $$
\end{theorem}
\texttt{... Beweis fehlt noch, wäre aber sehr praktisch ...}

\texttt{... hier noch eine Hausaufgabe für Thomas ... ;)}\todo{T: Eigentlich gew\"unscht ist $\ICC^s(\IHH,\IKK)\subset\ICC^\w(\IHH,\IKK)$, wobei nat\"urlich $\ICC^\w(\IHH,\IKK)=\{T\in\Gamma(\ILL(\IHH^\w,\IKK^\w)|$ $T\text{ strongly cont.}\}$ F\"ur $T\in\ICC^\s$ kann ich $T\in\ICC^\w$ bisher nur zeigen, falls alle $T(x)$ normal sind, was insbesondere $\IHH=\IKK$ voraussetzen muss...}
\begin{proposition}\label{stronggivesweak}
    Let $(\IHH,\ICC_\IHH)$ and $(\IKK,\ICC_\IKK)$ be continuous fields of Hilbert spaces. Then
    $$\elles(\IHH,\IKK)\subset \elles(\IHH^\w,\IKK^\w).$$
\end{proposition}
\begin{proof}
    The fibers of $\IHH^\w$ and $\IKK^\w$ are $\IHH_x$ and $\IKK_x$ equipped with their weak topologies, respectively, which we denote by $\IHH^\w_x$ and $\IKK^\w_x$, $x\in X$. As $\ILL(\IHH_x,\IKK_x)\subset \ILL(\IHH^\w_x,\IKK^\w_x)$ with continuous inclusion when the spaces of operators are equipped with the strong operator topology induced by the involved spaces, it trivially holds $\elles(\IHH,\IKK)\subset \elles(\IHH^\w,\IKK^\w).$ 
\end{proof}
\else
\section{A notion of strong convergence for closed operators on continuous fields of Banach spaces}
\label{sec:closedops}
In this section, we consider two continuous fields $(\IEE,\ICC_\IEE)$ and $(\IFF,\ICC_\IFF)$ of Banach spaces over the same topological space $X$. We will see that a weakening of the notion of \emph{strong graph limits}, see, e.g., \cite{resi}, Section VIII.7 has a natural extension to the setting of continuous fields and gives a notion of strong convergence that coincides with that defined in \ref{def:strongoperator convergence} in case the operators in question are bounded. At the same time, it extends what is known as \emph{strong resolvent convergence} or generalized convergence as considered in \cite{kato}, Chapter VIII, \S 1 .

To fix notation, denote by $\ICl(\IEE_0,\IFF_0)$ the set of closed operators between two Banach spaces $\IEE_0$ and $\IFF_0$. Of course an operator $T$ from $\IEE_0$ to $ \IFF_0$ with domain of definition $D(T)\subset \IEE_0$ is said to be \emph{closed}, when its graph
$$
G(T)=\set{(\zeta,T\zeta)}{\zeta\in D(T)}\subset \IEE_0\times\IFF_0
$$
is closed; when we consider it as a Banach space with the norm induced by $\IEE_0\times\IFF_0$ (built like the graph norm below so that it is a Hilbert norm in the case that $\IEE_0$ and $\IFF_0$ are Hilbert spaces) we write $\IGG(T)$.

The space $D(T)$ endowed with the \emph{graph norm}
$$
\|\cdot\|_T :=\left(\|\cdot\|^2+ \| T\cdot\|^2\right)^\frac{1}{2}
$$
is a Banach space denoted by $\IDD(T)$ and $T$ induces a bounded operator $T^\flat\in\ILL(\IDD(T),\IFF_0)$ by $T^\flat\zeta=T\zeta$.

\begin{definition} In the situation above:
\begin{itemize}
   \item [(1)] We denote by
   $$\ICl(\IEE,\IFF):=\set{T}{T\in\ICl(\IEE_x,\IFF_x)\mbox{  for some  }x:=\pi(T)}
   $$
   the total space of closed operators from $\IEE$ to $\IFF$.
   \item [(2)] A section $T\in\Gamma(\ICl(\IEE,\IFF))$ induces a field $\IDD(T)$ of Banach spaces and the section $T^\flat\in\Gamma(\ILL(\IDD(T),\IFF))$.
   \item [(3)] On $\IDD(T)$ we define
   $$
   \ICC(T):=\set{f\in\ICC_\IEE\cap\Gamma(\IDD(T))}{Tf\in\ICC_\IFF}
$$
   \item [(4)] We define a field of Banach spaces by
   $$
   \IGG(T):=\bigsqcup_{x\in X}\IGG(T_x)
   $$
   \item [(5)] A section $T\in\Gamma(\ICl(\IEE,\IFF))$ is said to be  \emph{G-continuous}, if for each $x\in X$ and any $\zeta\in D(T_x)$ there is $f\in \ICC(T)$ such that $f(x)=\zeta$.
   \item [(6)] A net $\snet{T}{\alpha}{A}$ in $\ICl(\IEE,\IFF)$ is said to \emph{G-converge to} $T_\infty\in \ICl(\IEE,\IFF)$ if $\pi(T_\alpha)\to \pi(T_\infty)$ and, for any $\zeta\in D(T_\infty)$ there exist $\zeta_\alpha\in D(T_\alpha)$ so that $\lim \zeta_\alpha=\zeta_\infty$ and $\lim T_\alpha\zeta_\alpha=T_\infty\zeta_\infty$. We then write $T_\alpha\xrightarrow{G} T_\infty$. 
   
\iffalse
\todo{Ali: Why is G-convergence well defined? It seems that the notion of G-convergence is not well defined, for the limit $\lim T_\alpha\zeta_\alpha$ depends on the choice of the net $\zeta_\alpha$ and leads to a multi-valued operator in general. I would propose to adapt the definition of Reed-Simon, p.293 (or  Ethier-Kurtz, Joergensen): We say that $T_\alpha$ G-converges to $T$ if $\pi(T_\alpha)\to \pi(T)$ and  the set
\[
 {\rm {ext-lim}} T_\alpha := \big\{ (\zeta,\xi),\ \pi(\zeta) = \pi(\xi)=\pi(T),\ \exists\,\zeta_\alpha\in D(T_\alpha),\ \zeta_\alpha\to\zeta\ \text{and}\ 
 T_\alpha\zeta_\alpha\to\xi \big\},
 \]
 coincides with $G(T)$. 
With this definition,  Prop.-(2), Coro.5.5-(2) and Prop.5.6-(2) are true under the assumption that $B_\alpha, T_\alpha^{-1}$ are eventually bounded. The latter assumption is very convenient for Part 6, as we are dealing with self-adjoint operators.\\
If you agree with this observation I will put the proofs in the file.}
\else
\end{itemize}
\end{definition}
\begin{remark}\label{rem:G-convergence}
    \begin{itemize}
    \item [\textrm{(1)}] Since $\ICC_{\IEE}$ is saturated, it is easily seen that every net $(T_\alpha)_{\alpha\in A}$ in $\ILL(\IEE,\IFF)$ which converges meta strongly to $T_\infty\in \ILL(\IEE,\IFF)$ also G-converges to $T_\infty$.
    \item [\textrm{(2)}] Our notion of G-convergence is weaker than the convergence in the strong graph limit sense, as introduced, e.g., in \cite{resi}, Section VIII.7. The latter can be rephrased as the Kuratowski convergence of the $G(T_\alpha)$ to $G(T_\infty)$, see, e.g., \cite{beer}. It includes the additional requirement that limit points of the graphs of the $T_\alpha$ belong to the graph of $T_\infty$. Our notion has the advantage that it fits perfectly with the framework of continuous fields, see Theorem \ref{thmGcont} below. Limits are by no means unique: if $T_\alpha\xrightarrow{G} T_\infty$ and $S$ is closed with $G(S)\subset G(T_\infty)$, then $T_\alpha\xrightarrow{G} S$ as well. However, this potential disadvantage of ambiguity does not bother us in the equivalences that we prove later due to the fact that $T_\infty$ or its inverse will be defined everywhere.
    \end{itemize}
\end{remark}
\begin{lemma}\label{lemmaclosed} Let $T\in\Gamma(\ICl(\IEE,\IFF))$.
 \begin{itemize}
   \item [(1)] $\ICC(T)$ is a locally uniformly closed continuous subspace for $\IDD(T)$.
   \item [(2)]  $\ICC(\IGG(T)):= \Gamma(\IGG(T))\cap (\ICC_\IEE\times \ICC_\IFF)$
   is a locally uniformly closed continuous subspace for $\IGG(T)$.
 \end{itemize}  
\end{lemma}
\begin{proof}
That $\ICC(T)$ and $\ICC(\IGG(T))$ are continuous for $\IDD(T)$ and $\IGG(T)$, respectively, is obvious, while the local uniform closedness follows immediately from the fact that $\ICC_\IEE$ and $\ICC_\IFF$ are locally uniformly closed.
\end{proof}
\begin{theorem}\label{thmGcont}
For $T\in\Gamma(\ICl(\IEE,\IFF))$ the following are equivalent:
\begin{itemize}
 \item [\textrm{(i)}] $T$ is G-continuous.
 \item [\textrm{(ii)}] $(\IDD(T),\ICC(T))$ is a continuous field of Banach spaces and $T^\flat\in\ICC^\s(\IDD(T),\IFF)$.
 \item [\textrm{(iii)}] $(\IGG(T),\ICC(\IGG(T)))$ is a continuous field of Banach spaces.
\end{itemize}
\end{theorem}
\begin{proof}
The G-continuity of $T$ implies that $\ICC(T)$ is saturated. By Lemma \ref{lemmaclosed}, (i) therefore implies that $(\IDD(T),\ICC(T))$ is a continuous field of Banach spaces and thus, $T^\flat\in\ICC^\s(\IDD(T),\IFF)$. Hence, (i) implies (ii). Again by Lemma \ref{lemmaclosed}, (iii) trivially follows from (ii). Finally, if (iii) holds, $\ICC(\IGG(T))$ is saturated for $\IGG(T)$ so that for any $\zeta\in D(T_x)$ there is $(f_1,f_2)\in \ICC(\IGG(T))= \Gamma(\IGG(T))\cap (\ICC_\IEE\times \ICC_\IFF)$ with $(f_1(x),f_2(x))=(\zeta, T_x \zeta)$. Thus, $T$ is G-continuous.
\end{proof}

\begin{proposition}\label{propbddperturbation}
\begin{itemize}
 \item [(1)] For all $T\in\Gamma(\ICl(\IEE,\IFF))$ and $B\in\ICC^\s(\IEE,\IFF)$  
 \iffalse
 \todo{Ali: hier soll  $B\in\ICC^\s(\IEE,\IFF)$ vielleicht; Thomas: ich denke, hier reicht $B\in\Gamma(\IEE,\IFF)$, oder?} 
 \else
 one has the equivalence
 $$T\mbox{  is G-continuous }\Longleftrightarrow T+B\mbox{  is G-continuous.}$$
 \item [(2)] Let $\snet{T}{\alpha}{A}$ be a net in $\ICl(\IEE,\IFF)$ and $T_\infty\in \ICl(\IEE,\IFF)$, let $\snet{B}{\alpha}{A}$ be a net in $\ILL(\IEE,\IFF)$, $B_\infty\in
 \ILL(\IEE,\IFF)$ with $\pi(T_\alpha)=\pi(B_\alpha)$ $(\alpha \in A)$ and $B_\alpha\xrightarrow{\ms}B_\infty$. Then one has the equivalence
 $$T_\alpha \xrightarrow{G}T_\infty\mbox{ }\Longleftrightarrow \mbox{ }T_\alpha+B_\alpha\xrightarrow{G}T_\infty+B_\infty$$
\end{itemize}
\end{proposition}
\begin{proof}
 (1) Assume that $T\in\Gamma(\ICl(\IEE,\IFF))$ is G-continuous  and $B\in\ICC^\s(\IEE,\IFF)$. Let $x\in X, \zeta\in D(T_x + B_x) = D(T_x)$, by assumption there is $f\in \ICC(T)$ with $f(x) = \zeta$. As $B\in\ICC^\s(\IEE,\IFF)$ we get $Bf\in\ICC_\IFF$,  $Tf + Bf\in\ICC_\IFF$ and hence $f\in \ICC(T + B)$ and satisfies $f(x)=\zeta$ so that $T+B$ is $G$-continuous.\\
 Taking into account the continuity of $\IFF\rightarrow\IFF, \eta\mapsto -\eta$ (see Remark \ref{rem:natural topology is an initial topology}), similar arguments lead to the converse implication.

 (2) follows from (1) via the net machine for the net $\overline{x}=\left(\pi(T_\alpha)\right)_{\alpha\in A\cup\{\infty\}}$.
\end{proof}
\begin{corollary}\label{G-Strong}
\begin{itemize}
 \item [(1)] Let $B\in\Gamma(\ILL(\IEE,\IFF))$ be locally bounded, i.e., for every $x_0\in X$ there is $U(x_0)\in\mathfrak{U}(x_0)$ with $\sup_{x\in U(x_0)}\|B(x)\|<\infty$
 \iffalse
 \todo{Thomas: Der Index ${\IEE_x\rightarrow\IFF_x}$ hier und in (2) ist sicher Notations-Overkill, den man in der endg\"ultigen Fassung beseitigen kann...}
 \else
 . Then one has the equivalence
 $$B\mbox{  is G-continuous }\Longleftrightarrow B\in\ICC^\s(\IEE,\IFF).$$
 \item [(2)] Let $\snet{B}{\alpha}{A}$ be a net in $\ILL(\IEE,\IFF)$ with $\sup_{\alpha\ge\alpha_0}\|B_\alpha\|<\infty$ for some $\alpha_0\in A$, $B_\infty\in
 \ILL(\IEE,\IFF)$. Then one has the equivalence
 $$B_\alpha\xrightarrow{G}B_\infty\mbox{ }\Longleftrightarrow \mbox{ }B_\alpha\xrightarrow{\ms}B_\infty .$$
 \end{itemize}
\end{corollary}
\begin{proof}
  For $B\in\Gamma(\ILL(\IEE,\IFF))$ we have $\IDD(B)=\IEE$ and $\ICC(B)=\set{f\in \ICC_{\IEE}}{B f\in\ICC_{\IFF}}$. Clearly, this implies the G-continuity of $B$ whenever $B\in\ICC^s(\IEE,\IFF)$. Conversely, if $B$ is G-continuous and $f\in\ICC_{\IEE}$, for arbitrary $x_0\in X$ there is $h_{x_0}\in\ICC(B)$ with $h_{x_0}(x_0)=f(x_0)$. Since $B h_{x_0}\in\ICC_{\IFF}$ and
  $$\|B(x)f(x)-B(x)h_{x_0}(x)\|\leq \sup_{x\in U(x_0)}\|B(x)\|\|f(x)-h_{x_0}(x)\|,$$
  it follows from $f-h_{x_0}\in \ICC_{\IEE}$ and the local uniform closedness of $\ICC_\IFF$ that $Bf\in\ICC_\IFF$ which implies $B\in\ICC^s(\IEE,\IFF)$. Thus, (1) holds. 

  To prove (2), we allude to Lemma \ref{strongconvergence} in order to see that $G$-convergence implies meta strong convergence. As observed in Remark \ref{rem:G-convergence}-(1), the converse implication is trivially true.
\end{proof}

\begin{proposition}\label{propGinvers} Let $S\in\Gamma(\ICl(\IEE,\IFF))$,  $\snet{T}{\alpha}{A}$ be a net in $\ICl(\IEE,\IFF)$ and $T_\infty\in \ICl(\IEE,\IFF)$. 
\begin{itemize}
 \item [(1)] If $S_x$ is injective for all $x\in X$, then $S^{-1}\in\Gamma(\ICl(\IFF,\IEE))$ and
 $$S\mbox{  is G-continuous }\Longleftrightarrow S^{-1}\mbox{  is G-continuous.}$$
 \item [(2)] If all $T_\alpha$, $T_\infty$ are injective, then one has the equivalence
 $$T_\alpha \xrightarrow{G}T_\infty\mbox{ }\Longleftrightarrow \mbox{ }T_\alpha^{-1}\xrightarrow{G}T_\infty^{-1}.
 $$
 \item [(3)] If $S_x:D(T_x)\to\IFF_x$ is bijective for all $x\in X$, and $S^{-1}$ is locally bounded, one has the equivalence 
 $$S\mbox{  is G-continuous }\Longleftrightarrow S^{-1}\in\ICC^\s(\IFF,\IEE)$$ 
 \iffalse
 \todo{Ali: ich denke hier soll $S^{-1}\in \ICC^s$ stehen}
 \else
 \item [(4)] If all $T_\alpha$, $T_\infty$ are bijective, meta strong convergence of $(T^{-1}_\alpha)_{\alpha\in A}$ to $T^{-1}_\infty$ implies $
 T_\alpha \xrightarrow{G}T_\infty$. Under the additional hypothesis $\sup_{\alpha\ge\alpha_0}\|T_\alpha^{-1}\|<\infty$ for some $\alpha_0\in A$, one has the equivalence
 $$
 T_\alpha \xrightarrow{G}T_\infty\mbox{ }\Longleftrightarrow \mbox{ }T_\alpha^{-1}\xrightarrow{\ms}T_\infty^{-1}.
 $$
\end{itemize}
\end{proposition}
\begin{proof}
(1) is obvious, (2) follows from (1) by the net machine, and (3) is a consequence of the Closed Graph Theorem and Corollary \ref{G-Strong}-(1). One part of (4) is a consequence of (2) and Remark \ref{rem:G-convergence}-(1) while the claim under the additional hypothesis holds true by Lemma \ref{strongconvergence}-(2).
\end{proof}

\section{Convergence of self-adjoint operators on the total space of a continuous field of Hilbert spaces}
\label{sec:self-adjointops}
Let again $(\IHH,\ICC)$ be a continuous field of Hilbert spaces over $X$. We abbreviate as usual $\ILL(\IHH)=\elle{H}{H}$ and regard it as the total space of a field of locally convex spaces and $\ICC^\s(\IHH) = \ICC^\s(\IHH,\IHH)$ the set of strongly continuous sections. We recall that according to Corollary \ref{Ls} if $X$ is a locally compact Hausdorff space or a normal space and $\IHH_x$ are separable for every $x$ then $(\ILL(\IHH),\ICC^\s(\IHH))$ is a continuous  field of locally convex spaces.

We first record a straightforward extension of the spectral consequences of strong convergence of bounded self-adjoint operators:
\begin{theorem}\label{thm:functionalcalculus} Assume $(T_\alpha)_{\alpha\in A}$ is a net of self-adjoint operators in $\ILL(\IHH)$ such that for some $\alpha_0\in A$ it holds $\sup_{\alpha\ge \alpha_0} \|T_\alpha\|<\infty$, and let $T_\infty\in \ILL(\IHH)$.\\ 
\emph{(1)} One has $T_\alpha\xrightarrow{\ms} T_\infty$, if and only if for all $\phi\in C(\C)$ one has $\phi(T_\alpha) \xrightarrow{\ms} \phi(T_\infty)$.\\
\emph{(2)} If $T_\alpha\xrightarrow{\ms} T_\infty$, then $\sigma(T_\infty)\subset \lim_{\alpha} \sigma(T_\alpha)$, in the sense that for all $\lambda\in \sigma(T_\infty)$ there exists a net $(\lambda_\alpha)_{\alpha\in A}$ in $\R$ with $\lambda_\alpha\in \sigma(T_\alpha)$ for all $\alpha\in I$, such that $\lambda_\alpha\to \lambda$.
\end{theorem}

\begin{proof}
(1) The 'if' is evident and the 'only if' implication can be proven by a standard Stone-Weierstrass argument: due to the fact that strong convergence does not pass to adjoints, a little caution is required. We can assume without restriction that $\sup_{\alpha\in A} \|T_\alpha\|<\infty$ and so the spectra of the operators involved are contained in a compact interval $[a,b]\subset\R$. As $\phi(T_\alpha), \phi(T_\infty)$ depend solely on the values of $\phi$ on $[a,b]$, it suffices to consider $\phi\in C([a,b];\C)$.  By linearity it suffices to prove the asserted convergence for  $\phi\in C([a,b];\R)$. To this end, consider 
$$
B:=\set{\phi\in C([a,b];\R)}{\phi(T_\alpha) \xrightarrow{\ms} \phi(T_\infty)}.
$$
Obviously, $B$ is a norm-closed subalgebra of $C([a,b];\R)$. By the assumed convergence $T_\alpha\xrightarrow{\ms} T_\infty$, it contains the identity and therefore separates the points of $[a,b]$. Moreover, the constant $1\in B$, so that $B$ vanishes nowhere. The Stone-Weierstrass theorem gives that $B=C([a,b];\R)$, proving the asserted convergence.

(2) Functional calculus gives that $\lambda\in\sigma(T_\infty)$ if and only if $\phi(T_\infty)\not= 0$ for every $\phi\in C(\C)$ with $\phi(\lambda)\not=0$. Together with (1) this gives the claim.
\end{proof}

We will be concerned with strong continuity and convergence of self-adjoint operators over $\IHH$. We can build on the notions developed in Sections  \ref{sec:opsbetweencontfields} and \ref{sec:closedops} and will also use the analysis provided in Section \ref{sec:weakandstrong}. To ease notation we introduce the total space
$$
\Isad(\IHH):=\set{T}{T\mbox{  self-adjoint in } \IHH_x\mbox{  for some } x:=\pi(T)}
$$

\begin{corollary}\label{reso}  Assume $(T_\alpha)_{\alpha\in A}$ is a net of self-adjoint operators in $\Isad(\IHH)$, and that $T_\infty\in \Isad(\IHH)$. Then the following statements are equivalent:
\begin{itemize}
\item[(i)] There exists 
$$
\lambda \in  \bigcap_{\alpha } \varrho(T_\alpha)\cap \varrho(T_\infty)
$$
with 
$$(\lambda-T_\alpha)^{-1}\xrightarrow{\ms} (\lambda-T_\infty)^{-1}.$$
\item[(ii)] For all 
$$
\lambda \in  \bigcap_{\alpha } \varrho(T_\alpha)\cap \varrho(T_\infty)
$$
one has
$$(\lambda-T_\alpha)^{-1}\xrightarrow{\ms} (\lambda-T_\infty)^{-1}.$$
\item[(iii)] For all $\phi\in  C_0(\C)$ one has $\phi(T_\alpha) \xrightarrow{\ms} \phi(T_\infty)$.
\end{itemize}
\end{corollary}

\begin{proof}
(i)$\Longrightarrow$(iii): We use the trivial implication in Proposition \ref{propGinvers} (4) to conclude that $(\lambda-T_\alpha)^{-1}\xrightarrow{G} (\lambda-T_\infty)^{-1}$. By Proposition \ref{propGinvers} (2) we get  $\lambda-T_\alpha\xrightarrow{G} \lambda-T_\infty$. An appeal to Proposition \ref{propbddperturbation} gives that $\lambda'-T_\alpha\xrightarrow{G} \lambda'-T_\infty$ for every $\lambda'\in \C$. Specializing to $\lambda'=ci$ for $c\in\R\setminus\{0\}$ and using Proposition \ref{propGinvers} (4) again, we get that $(ci-T_\alpha)^{-1}\xrightarrow{\ms} (ci-T_\infty)^{-1}$. 

We are now ready to do a Stone-Weierstrass gavotte as in the proof of the preceding theorem:
First note that in (iii) we actually need to deal with $\phi\in C_0(\R;\C)$ as the spectra of the operators involved are real. By linearity it suffices to prove the asserted convergence for  $\phi\in C_0(\R;\R)$. To this end, consider 
$$
B:=\set{\phi\in C_0(\R;\R)}{\phi(T_\alpha) \xrightarrow{\ms} \phi(T_\infty)}.
$$
Obviously, $B$ is a norm-closed subalgebra of $C_0(\R;\R)$. It contains the functions 
$$f_c(\cdot):=(\cdot-ci)^{-1}+(\cdot+ci)^{-1}, c\in\R\setminus\{0\}$$ 
by our discussion above. Since these functions separate the points of $\R$, so does $B$. Moreover, $g(\cdot):=(\cdot -i)^{-2} + (\cdot +i)^{-2}\in B$, so that $B$ vanishes nowhere. The Stone-Weierstrass theorem gives that $B=C_0(\R;\R)$, proving the asserted implication.

(iii)$\Longrightarrow$(ii) and (ii)$\Longrightarrow$(i) are clear.
\end{proof}

\begin{definition} In the situation of Corollary \ref{reso}, we say that $(T_\alpha)_{\alpha\in A}$ \emph{converges to $T_\infty$ in the strong resolvent sense}. We then write $T_\alpha \xrightarrow{\srs} T_\infty$.
\end{definition}
From this definition and Proposition \ref{propGinvers} we immediately get:
\begin{corollary}\label{reso=G}  Assume $(T_\alpha)_{\alpha\in A}$ is a net of self-adjoint operators in $\Isad(\IHH)$, and that $T_\infty\in \Isad(\IHH)$. Then 
$$
T_\alpha \xrightarrow{\srs} T_\infty\quad\Longleftrightarrow\quad T_\alpha \xrightarrow{G} T_\infty
$$
\end{corollary}

Theorem \ref{thm:functionalcalculus} also immediately implies:

\begin{corollary}
 If $T_\alpha \xrightarrow{\srs} T_\infty$, then $\sigma(T_\infty)\subset \lim_\alpha \sigma(T_\alpha)$.
\end{corollary}

\section*{Mosco-Convergence revisited}\vspace{.5cm}

In this subsection, we state and prove a useful characterization of strong resolvent convergence for self-adjoint, nonnegative operators on a continuous field of Hilbert spaces, in terms of Mosco convergence of their quadratic forms. Typically, these forms are technically much easier to handle and this is why such a result is interesting. Similar results have been proved in various special cases before, cf. \cite{kuwae}.  

We now specialize further to nonnegative self-adjoint operators and define

$$
\Isaplus(\IHH):=\set{T}{T\mbox{  self-adjoint in } \IHH_x\mbox{  for some } \pi(T):=x, T\ge 0}
$$

For any  $T\in \Isaplus(\IHH)$ we denote by $\form{t}$ the associated closed form with domain $D(\form{t}):=D(T^\frac12)$ and inner product $\langle\cdot,\cdot\rangle_\form{t}:= \form{t}[\cdot,\cdot] + \langle\cdot,\cdot\rangle$, as usual we write 
\[
\form{t}[u] := \form{t}[u,u]\mbox{  for  }u\in D(\form{t}), 
\]
see \cite{kato}, Chapter VI, \S 2 for the background on forms and their associated operators. The important point to recall is that $D(\form{t})$ equipped with the form norm is a Hilbert space that we denote by $\IDD(\form{t})$. In view of the above mentioned equivalence we have a one-to-one correspondence between  $T\in\Gamma(\Isaplus(\IHH))$ and $\form{t}\in\Gamma\left(\IFF(\IHH)\right)$, where
$$\IFF(\IHH):=\set{\form{t}}{\form{t}\mbox{ is a densely defined, closed, nonnegative form in }\IHH_x\mbox{ for some }x=:\pi(\form{t})} .
$$
We define
\[
\IDD(\form{t}) : = \bigsqcup_{x\in X} \IDD(\form{t}_x),
\]
which consequently is a field of Hilbert spaces with the natural choice of $\langle\cdot,\cdot\rangle_{\form{t}_x}$ on the fibers, where it is understood that $\form{t}_x$ is associated with $T_x$ in the above sense. There is a natural candidate for a continuous subspace on $\IDD(\form{t})$, namely
\[
\mathscr{C}(\form{t}):= \set{f\in\mathscr{C}}{ f(x)\in D(\form{t}_x) \text{  for all}\ x\ \text{and}\ 
\form{t}[f]:=\form{t}_{(\cdot)}[f(\cdot)] \in C(X)} .
\]
    It is clear that  $\ICC(\form{t})$ is locally uniformly closed.

\begin{lemma}\label{lem:Q} In the above situation if $\ICC(T^\frac12)$ is saturated, then
$$
\ICC(\form{t})=\ICC(T^\frac12) ,
$$
and $(\IDD(\form{t}),\ICC(\form{t}))$ is a continuous field.
\end{lemma}
\begin{proof}
We only have to show that for $f\in\ICC(\form{t)}$ it follows that $T^\frac12 f\in\ICC$. Pick such an $f$ and $x_\infty\in X$. For any net $\snet{x}{\alpha}{A}$ that converges to $x_\infty$, Proposition \ref{prop:bddgivesmw} applied to  $(T^\frac12 f(x_\alpha))_{\alpha\in A}$ provides us with a  subnet $\varphi:B\to A$  so that $(T^\frac12 f(x_\beta))_{\beta\in B}$ is meta weakly convergent in $\IHH$. Since $\| f\|_\form{t}$ is continuous, an appeal to Proposition \ref{prop:mwandnormgivesstrong} gives that $(f(x_\beta))_{\beta\in B}$ and $(T^\frac12 (x_\beta)f(x_\beta))_{\beta\in B}$ are strongly convergent to $f(x_\infty)$ and $\xi$, respectively. 
We have to check that $\xi=T(x_\infty)^\frac12 f(x_\infty)$. To do so, pick $\zeta\in D(\form{t}(x_\infty))$ and $g\in \ICC(T^\frac12)$ with $g(x_\infty)=\zeta$ which is possible by assumption. By polarization we get that $\langle f,g\rangle_\form{t}$ is continuous and therefore
\begin{align*}
    \langle T(x_\infty)^\frac12 f(x_\infty), T(x_\infty)^\frac12 \zeta\rangle+\langle f(x_\infty),\zeta\rangle&=\langle f(x_\infty),\zeta\rangle_{\form{t}(x_\infty)}\\
    &=\lim_{\beta\in B}\langle f(x_\beta),g(x_\beta)\rangle_{\form{t}(x_\beta)}\\
    &=\lim_{\beta\in B}\left(\langle T(x_\beta)^\frac12 f(x_\beta),T(x_\beta)^\frac12 g(x_\beta)\rangle+\langle f(x_\beta),g(x_\beta)\rangle\right)\\
    &= \langle \xi , T(x_\infty)^\frac12 \zeta\rangle+\langle f(x_\infty),\zeta\rangle ,
\end{align*}
and since $\zeta$ was arbitrary in $D(T(x_\infty)^\frac12)$ this yields $\xi=T(x_\infty)^\frac12 f(x_\infty)$.
To sum up, we have shown that any net $\snet{x}{\alpha}{A}$ that converges to $x_\infty$ has a subnet  so that $\lim_{\beta\in B}T(x_\beta)^\frac12 f(x_\beta)=T(x_\infty)^\frac12 f(x_\infty)$ which implies the asserted continuity.
\end{proof}

The following fundamental result generalizes a number of results \cite{attouch,mosco,kuwae,kato,kim, chen}, on the equivalence of Mosco convergence and strong resolvent convergence. It also displays quite nicely the power of the concept of continuous fields, even in the case in which all operators live in one Hilbert space: namely, two characterizations of continuous dependence on a parameter of unbounded operators which are intrinsically in terms of continuous fields of Hilbert spaces.

\begin{theorem} 
\label{Cfield-Mcontinuity}
For $T\in\Gamma(\Isaplus(\IHH))$ and $\form{t}$ as above, the following assertions are equivalent
\begin{itemize}
\item[(i)] $\big(T + \lambda\big)^{-1}\in \ICC^\s(\IHH)$ for some $\lambda\in (0,\infty)$.
\item[(ii)] $(\IDD(T),\ICC(T))$ defines a continuous field of Hilbert spaces.
\item[(iii)] $(\IDD(\form{t}),\ICC(\form{t}))$ defines a continuous field of Hilbert spaces that is compatible with $\IHH^\w$ in the following sense:\\
$\rm{(\widehat{M1})}$ $\ICC(\form{t})^\w\subset\ICC^\w$.
\item[(iv)] $(\IDD(\form{t}),\ICC(\form{t}))$ defines a continuous field of Hilbert spaces that is compatible with $\IHH^\w$ in the following sense:\\
$\rm{(\overline{M1})}$ Every bounded $u\in \ICC(\form{t})^\w$ lies in $\ICC^\w$.

\item[(v)] $\form{t}$ is \emph{Mosco-continuous} in the sense that\\
$\rm{(\widetilde{M1})}$ If $u\in \ICC^\w$, then $(\form{t}+1)[u]$ is lower semicontinuous.\\ 
$\rm{(\widetilde{M2})}$ For any $\zeta\in \IDD(\form{t})$ there is $f\in \ICC(\form{t})$ such that $f(\pi(\zeta)) = \zeta$. 
\end{itemize}
\end{theorem}
We adopt the common convention that $\form{s}[u]:=\infty$, if $u\notin D(\form{s})$ which explains statement $\rm{(\widetilde{M1})}$. 

As in the proof by Mosco from \cite{mosco}, we use the following elementary but useful variational characterization of the resolvent. It is usually phrased for real Hilbert spaces, so we include it for the readers convenience:
\begin{proposition}\label{prop:variational}
For $\form{s}\in \IFF(\IHH)$, $S$ the associated self-adjoint operator and 
$\lambda\in\rho(S)\cap [0,\infty)$, the vector $u \in \IHH_{\pi(\form{s})}$, $(S+\lambda)^{-1} u$ is the unique minimizer of the functional
\begin{equation}\label{eq:variational}
D(\form{s})\ni v\mapsto (\form{s}+\lambda)[v]-2\Re\langle u,v\rangle . 
\end{equation}   
\end{proposition}
\begin{proof}
    Evidently, 
    \begin{align*}
        0&\le (\form{s}+\lambda)[(S+\lambda)^{-1}u-v]\\
        &= \langle (S+\lambda)^{-1}u,u\rangle +(\form{s}+\lambda)[v]-2\Re\langle u,v\rangle 
    \end{align*}
    has its unique minimum for $v=(S+\lambda)^{-1}u$, so the claim follows by forgetting the constant first term in the latter sum.
\end{proof}
\begin{proof}[Proof of Theorem \ref{Cfield-Mcontinuity}]

(i)$\Leftrightarrow$ (ii): is even true for self-adjoint operators, see Theorem \ref{thmGcont}, Proposition \ref{propGinvers} and Proposition \ref{propbddperturbation} above.

(i) $\Rightarrow$ (iii): follows with functional calculus. Since $(T+1)^{-1}\in  \ICC^\s(\IHH)$, so is $(T+1)^{-\frac12}\in  \ICC^\s(\IHH)$ by Theorem \ref{thm:functionalcalculus}. This implies that $\ICC(T^\frac12)$ is saturated and an appeal to Lemma \ref{lem:Q} gives the first part of the assertion.

To prove $\rm{(\widehat{M1})}$, let $u\in \ICC(\form{t})^\w$ and $g\in\ICC$. It follows that  $h:=(T+1)^{-1}g\in\ICC(T)$ and $h\in\ICC(T^\frac12)\subset \ICC(\form{t})$.

This gives
\begin{align*}
    \langle u,g\rangle &=\langle u,(T+1)h\rangle =\langle u,h\rangle_\form{t} \in C(X),
\end{align*}
which was to be proven.

(iv) $\Rightarrow$ (v): $\rm{(\widetilde{M2})}$ is just the fact that $\ICC(\form{t})$ is saturated.

To prove $\rm{(\widetilde{M1})}$, let $u\in\ICC^\w$; we argue with nets, i.e., we have to show that for every net $\snet{x}{\alpha}{A}$ converging to $x_\infty$ such that $\lim (\form{t}_{x_\alpha}+1)[u(x_\alpha)]=:\tau_\infty<\infty$ exists, it follows that $(\form{t}_{x_\infty}+1)[u(x_\infty)]\le \tau_\infty$.
Using the net machine, we can work with $\IHH_{\overline{x}}$, $\IDD(\form{t})_{\overline{x}}$ and the corresponding continuous subspaces. Observe that due to the boundedness, Proposition \ref{prop:bddgivesmw} gives a subnet $B\xrightarrow{\varphi}A$ so that $(u(x_{\varphi(\beta)}))_{\beta\in B}$ converges in $\ICC(\form{t})^\w$. By $\rm{(\overline{M1})}$, its limit is $u(x_\infty)$ and the lower semicontinuity of the norm w.r.t weak convergence in $\IDD(\form{t})^\w$, Corollary \ref{cor:lscnorm}, implies that 
$\|u(x_\infty)\|^2_{\form{t}(x_\infty)}= (\form{t}_{x_\infty}+1)[u(x_\infty)]\le \tau_\infty$, as was to be shown.

(v) $\Rightarrow$ (i):  We start with $u\in\ICC$ and aim to show that $(T+1)^{-1}u\in\ICC$. To this end, pick $g\in\ICC(\form{t})$; since $g\in\ICC$, it follows that
$$
    \langle (T+1)^{-1}u,g\rangle_\form{t}
=\langle u,g\rangle \in C(X), 
$$
meaning that $(T+1)^{-1}u\in\ICC(\form{t})^\w$. We will now first prove that $(T+1)^{-1}u\in\ICC^\w$; to this end, we use the variational characterization of the resolvent from
Proposition \ref{prop:variational} above. Again, we use nets; we have to show that for every net $\snet{x}{\alpha}{A}$ converging to $x_\infty$, 
$$
(T_{x_\alpha}+1)^{-1}u(x_\alpha)\to (T_{x_\infty}+1)^{-1}u(x_\infty)\mbox{  in  }\IHH^\w .
$$
This, in turn, amounts to verifying that for every subnet of $\snet{x}{\alpha}{A}$ (which we denote by the same symbol to ease notation) there is a subnet $B\xrightarrow{\varphi}A$ so that
$$
(T_{x_{\varphi(\beta)}}+1)^{-1}u(x_{\varphi(\beta)})\to (T_{x_\infty}+1)^{-1}u(x_\infty)\mbox{  in  }\IHH^\w .
$$
We abbreviate $T_\beta:=T_{x_{\varphi(\beta)}}$ and use indices for forms in the same way, as well as $u_\beta:= u(x_{\varphi(\beta)})$ and $T_\infty:=T_{x_\infty}$ and $u_\infty:= u(x_\infty)$. By the usual boundedness argument, we can pick the subnet in such a way that
\begin{equation}\label{eq:1proofCfield-Mcontinuity}
(T_\beta+1)^{-1}u_\beta\to w_\infty\mbox{  in  }\IHH^\w 
\end{equation}
for some $w_\infty\in \IHH$. By \eqref{eq:variational}, above, for arbitrary $\beta\in B$,
\begin{equation}\label{eq:2proofCfield-Mcontinuity}
(\form{t}_\beta+1)[(T_\beta+1)^{-1}u_\beta]-2\Re\langle u_\beta, (T_\beta+1)^{-1}u_\beta\rangle\le (\form{t}_\beta+1)[v_\beta]-2\Re\langle u_\beta, v_\beta\rangle \mbox{  for all }v_\beta .   
\end{equation}
We now pick $v_\infty\in D(\form{t}_\infty)$ and $v\in \ICC(\form{t})$ with $v(x_\infty)=v_\infty$ which is possible, since $\ICC(\form{t})$ is saturated and put   
$v_\beta:= v(x_{\varphi(\beta)})$. Taking the $\liminf$ on both sides of \eqref{eq:2proofCfield-Mcontinuity}, taking into account property $\rm{(\widetilde{M1})}$ and \eqref{eq:1proofCfield-Mcontinuity}, we get
\begin{equation}\label{eq:3proofCfield-Mcontinuity}
(\form{t}_\infty+1)[w_\infty]-2\Re\langle u_\infty, w_\infty\rangle\le (\form{t}_\infty+1)[v_\infty]-2\Re\langle u_\infty, v_\infty\rangle .
\end{equation}
Since $v_\infty\in D(\form{t}_\infty)$ was arbitrary, another appeal to Proposition \ref{prop:variational} shows that $w_\infty= (T_\infty+1)^{-1}u_\infty$, concluding the proof
of $(T+1)^{-1}u\in\ICC^\w$. This now gives that
$$
\langle (T+1)^{-1}u,(T+1)^{-1}u\rangle_\form{t}=\langle (T+1)^{-1}u,u\rangle \in C(X)
$$
and since  $(T+1)^{-1}u\in\ICC(\form{t})^\w$, we get $(T+1)^{-1}u\in\ICC(\form{t})$ by using Proposition \ref{prop:mwandnormgivesstrong}. Since $\ICC(\form{t})\subset\ICC$ this finishes the proof
of (iv) $\Rightarrow$ (i).
\end{proof}

The names $\rm{(\widetilde{M1})}$ and  $\rm{(\widetilde{M2})}$ are chosen in accordance with the usual terminology regarding Mosco-convergence. 
The above characterization of continuity of sections has an immediate counterpart in terms of convergence of nets in the total space $\IFF(\IHH)$.

\begin{definition} Let  $(\form{t}_\alpha)_{\alpha\in A}$ be a net in $\IFF(\IHH)$,
$\form{t}_\infty\in\IFF(\IHH)$ with associated self-adjoint operators $T_\alpha, \alpha\in A$,  $T_\infty$. 

We say that  $(\form{t}_\alpha)_{\alpha\in A}$ in $\IFF(\IHH)$ \emph{converges to $\form{t}_\infty$ in the sense of Mosco}, if $\pi(\form{t}_\alpha)\to \pi(\form{t}_\infty)$ and the following  properties are satisfied:\\
(M1) For $\snet{\zeta}{\alpha}{A}$ with $\pi(\zeta_\alpha)=\pi(\form{t}_\alpha)$ it follows that
$$
\zeta_\alpha\xrightarrow{\w}\zeta_\infty
\Longrightarrow
(\form{t}_\infty+1)[\zeta_\infty]\le \liminf (\form{t}_\alpha+1)[\zeta_\alpha] .$$
(M2) For all $\eta_\infty\in D(\form{t}_\infty)$ there exist $\eta_\alpha\in D(\form{t}_\alpha)$ for all $\alpha\in A$ 
with
    $$\eta_\alpha\xrightarrow{\s}\eta_\infty\mbox{  and  }\form{t}_\alpha[\eta_\alpha]\to\form{t}_\infty[\eta_\infty] .$$
We then write $\form{t}_\alpha\xrightarrow{\Mos} \form{t}_\infty$.
\end{definition}

The following result is an immediate consequence of the previous main Theorem \ref{Cfield-Mcontinuity} by virtue of the net machine 

\begin{corollary}\label{mosco} Let  $(\form{t}_\alpha)_{\alpha\in A}$ be a net in $\IFF(\IHH)$,
$\form{t}_\infty\in\IFF(\IHH)$, $T_\alpha$ the self-adjoint  operator associated with $\form{t}_\alpha$ for $\alpha\in A\cup\{\infty\}$. Then: 
\begin{itemize}
    \item[(1)] 
 $(T_\alpha+1)^{-1}\xrightarrow{\ms} (T_\infty+1)^{-1}\Longleftrightarrow \form{t}_\alpha\xrightarrow{\Mos} \form{t}_\infty$.
\item[(2)] For Mosco convergence \emph{(M1)} can be replaced by each of the following: 
\begin{itemize}
    \item [(M1*)] If $u_\alpha\in D(\form{t}_\alpha)$ is bounded w.r.t $\|\cdot\|_{\form{t}(\alpha)}$ and $u_\alpha\xrightarrow{\form{t}-\w} u_\infty$ in the sense that for all $v_\alpha\in D(\form{t}_\alpha)$ with $v_\alpha\xrightarrow{\s} v_\infty$ and $\form{t}_\alpha[v_\alpha]\to\form{t}_\infty[v_\infty]$, it follows that   
    $$
    \langle u_\alpha, v_\alpha\rangle_{\form{t}(\alpha)}\to \langle u_\infty, v_\infty\rangle_{\form{t}(\infty)}.
    $$
    Then $ u_\alpha\xrightarrow{\w}u_\infty$.
    \item[(M1**)]  If $u_\alpha\in D(\form{t}_\alpha)$ is bounded w.r.t $\|\cdot\|_{\form{t}(\alpha)}$ and $u_\alpha\xrightarrow{\form{t}-\w} u_\infty$ as well as $ u_\alpha\xrightarrow{\w}w_\infty$, then $u_\infty=w_\infty$.
\end{itemize}
\item[(3)]  For Mosco convergence, \emph{(M2)} can be replaced by
\begin{itemize}
    \item [(M2*)] For all $\eta_\infty$ in a dense subspace $D\subset D(\form{t}_\infty)$ there exist $\eta_\alpha\in D(\form{t}_\alpha)$ for all $\alpha\in A$ 
with
    $$\eta_\alpha\xrightarrow{\s}\eta_\infty\mbox{  and  }\form{t}_\alpha[\eta_\alpha]\to\form{t}_\infty[\eta_\infty] .$$
\end{itemize}
\end{itemize}
\end{corollary}
\begin{proof}
We use the net machine: consider  the corresponding continuous field $(\IHH_{\overline{x}},\ICC_{\overline{x}})$; then the net $(\form{t}_\alpha;\alpha\in A\cup\{\infty\})$ can be viewed as a section in $\Gamma(\IFF(\IHH_{\overline{x}}))$. From the definition we get that
\begin{align*}
\ICC_{\overline{x}}&=\set{v}{v_\alpha\xrightarrow{\s} v_\infty}\\
\ICC(\form{t})&=\set{v}{v_\alpha\xrightarrow{\s} v_\infty\mbox{  and  }\form{t}_\alpha[v_\alpha]\to\form{t}_\infty[v_\infty]}\\ 
\ICC(\form{t})^\w &=\set{v}{v_\alpha\xrightarrow{\form{t}-\w} v_\infty}\\ 
\end{align*}

To prove (1) it remains to observe that (M1) is just a reformulation of $\rm{(\widetilde{M1})}$ and (M2)  is  a reformulation of $\rm{(\widetilde{M2})}$ and apply Theorem \ref{Cfield-Mcontinuity}. 

(2) (M1*)   is just a reformulation of $\rm{(\widehat{M1})}$. Clearly, (M1*)$\Longrightarrow$(M1**). For the converse assume that $u_\alpha\xrightarrow{\form{t}-\w}u_\infty$. To show that $u_\alpha\xrightarrow{\w}u_\infty$, we have to prove that every subnet of the latter has a subnet that converges to $u_\infty$. To ease notation, let us denote the subnet by $u_\alpha$. By the  boundedness we assumed, an appeal to Proposition \ref{prop:bddgivesmw} gives a subnet $\varphi:B\to A$ so that $u_{\varphi(\beta)}\xrightarrow{\w}w_\infty$ for some $w_\infty$. By (M1**), it follows that $u_\infty=w_\infty$, giving the desired convergence.   
With this being said, the assertion of (2) again follows from Theorem \ref{Cfield-Mcontinuity}.

(3) Clearly, (M2) implies (M2*). For the converse recall that $\ICC(\form{t})$ is locally uniformly closed, so that (M2*) is sufficient to see that it is saturated, which is (M2).

\iffalse
(i)$\Longrightarrow$(ii): To prove (M1) we may assume by passing to a subnet that 
$$
C:= \sup_{\alpha\in A}\form{t}_\alpha[\zeta_\alpha]<\infty
$$
and we have to check that $\form{t}_\infty[\zeta_\infty]\le C$. We use Mosco's original method and introduce 
With the respective approximations for the $\form{t}_\alpha$ and strong resolvent convergence we get
$$
\form{t}_\infty^{(\beta)}[\zeta_\infty]=\lim_\alpha\form{t}_\alpha^{(\beta)}[\zeta_\alpha]\le C,
$$
giving the claim. For the proof of (M2), we take $\xi_\infty:= T_\infty^{\frac12}\eta_\infty$. Pick a net s.t. $\xi_\alpha\xrightarrow{\s}\xi_\infty$. Then $\eta_\alpha:=T_\alpha^{-\frac12}\xi_\alpha$ does the job, since $T_\alpha^{-\frac12}\xrightarrow{\ms}T_\infty^{-\frac12}$ by the assumption and Theorem \ref{thm:functionalcalculus}.

(ii)$\Longrightarrow$(i) is deduced as (iv)$\Longrightarrow$(i) in the proof of Theorem \ref{Cfield-Mcontinuity}.
\else
\end{proof}
\begin{remark} 
\begin{enumerate}[(1)]
    \item If the forms involved are induced by bounded self-adjoint operators $B_\alpha$ in the sense that $\form{t}_\alpha[\cdot,\cdot]=\langle B_\alpha \cdot, \cdot\rangle$ for all $\alpha\in A\cup\{\infty\}$, then  $\form{t}_\alpha\xrightarrow{\Mos} \form{t}_\infty$ if and only if $B_\alpha\xrightarrow{\m\s}B_\infty$.
    \item Define
$$
\form{t}_\alpha^{(\beta)}[u,v]:= \beta\langle u-\beta(T_\alpha +\beta)^{-1}u,v\rangle\mbox{  for   }\alpha\in A\cup\{\infty\},\beta>0,
$$
the Yosida-Moreau approximation. It is well-known and easy to check that
$$
\form{t}_\alpha[u] = \lim\form{t}_\alpha^{(\beta)}[u]=\sup_{\beta>0}\form{t}_\alpha^{(\beta)}[u] \mbox{  for all  }u\in \IHH_{\pi(\form{t}_\alpha)} .
$$
\item Consequently, the preceding Corollary shows that  $\form{t}_\alpha\xrightarrow{\Mos} \form{t}_\infty$
if and only if  $\form{t}^{(\beta)}_\alpha\xrightarrow{\Mos} \form{t}^{(\beta)}_\infty$ for some, respectively all $\beta>0$.
\end{enumerate}
\end{remark}
We next combine Corollary \ref{mosco} with sequential continuity to obtain:
\begin{corollary}\label{cor:moscoseq} Let $X$ be first countable and $\IHH$ be a continuous field of Hilbert spaces over $X$. For $\form{t}\in\Gamma(\IFF(\IHH))$ and $(T+1)^{-1}\in\Gamma(\ILL(\IHH))$ the associated section as above, the following are equivalent: 
    \begin{itemize}
        \item [(i)] $(T+1)^{-1}\in\ICC^\s(\IHH)$.
        \item[(ii)] For every $x_\infty\in X$ and every sequence $\snet{x}{n}{\N}$ in $X$ so that $x_n\to x_\infty$ one has
        $$
        \form{t}_{x_n} \xrightarrow{\Mos}\form{t}_{x_\infty}.
        $$
        \item [(iii)] $(\IDD(\form{t}),\ICC(\form{t}))$ satisfies:\\
        \begin{itemize} \item[(M1')] For any bounded $f\in \ICC(\form{t})^\w$ and every sequence $\snet{x}{n}{\N}$ in $X$ so that $x_n\to x_\infty$ one has $$f(x_n)\xrightarrow{\ICC^\w} f(x_\infty).$$
            \item [(M2')] For every $x\in X$, $\set{f(x)}{f\in \ICC(\form{t})}$ is dense in $\IDD(\form{t}_x)$.
        \end{itemize}
    \end{itemize}
\end{corollary}
\iffalse\begin{proof}
(i)$\Longrightarrow$ (iv) follows from Theorem \ref{Cfield-Mcontinuity}.\\
(iv)$\Longrightarrow$ (i) Since $\sup_{x\in X} \|( T_{\form{t}_{x}}  + 1)^{-1}\|\leq 1$ and $X$ is first countable, we infer from Corollary \ref{cor:seqcontinuity} that (i) $\Longleftrightarrow ( T_{\form{t}_{x_n}}  + 1)^{-1}\xrightarrow{\ms} (T_{\form{t}_{x_\infty}}  + 1)^{-1}$ for any $x_\infty\in X$ and any $x_n\in X$ converging to $x_\infty$.\\
To prove the latter property we shall use Lemma \ref{strongconvergence}-(3). Without loss of generality, we assume that $t_{x_\infty}$ is densely defined so that   the subspace $D$ is dense in $\IHH_{x_\infty}$. Pick $\zeta_\infty\in D$ and let $\zeta_n$ be the sequence from (M2'). We have to show $(T_{\form{t}_{x_n}} + 1)^{-1}\zeta_n \xrightarrow{\s} (T_{\form{t}_{x_\infty}}  + 1)^{-1}\zeta_\infty$. Now the proof of the latter follows as the proof of step (iv)$\Longrightarrow$ (i) from Theorem \ref{Cfield-Mcontinuity}, however, by using in the step from (\ref{eq:2proofCfield-Mcontinuity}) to (\ref{eq:3proofCfield-Mcontinuity}), the assumption (M2') instead of the fact that $\ICC(\form{t})$ is saturated.\\
(i)$\Longleftrightarrow$ (iii) Using Corollary \ref{cor:seqcontinuity} once again,  we obtain from the very definition of $\ICC^\s(\IHH)$ and the saturation property for $\IHH$ in conjunction Theorem \ref{mosco}, for any $f\in\Gamma(\IHH)$ such that $f(x_n) \xrightarrow{\s} f(x_\infty)$,
\begin{align*}
(T_{\form{t}} + 1)^{-1}\in\ICC^\s(\IHH)\Leftrightarrow (T_{\form{t}_{x_n}} + 1)^{-1}f(x_n)\xrightarrow{\s} (T_{\form{t}_{x_\infty}} + 1)^{-1}f(x_\infty)
\Leftrightarrow \form{t}_{x_n} \xrightarrow{\Mos}\form{t}_{x_\infty}.
\end{align*}
The equivalence (i)$\Longleftrightarrow$ (ii) follows by similar arguments in conjunction with Theorem \ref{Cfield-Mcontinuity}.
\end{proof}
\else
We end this section with an instructive example of how (M1) can fail:
\begin{example}
    Let $X=\N\cup\{\infty\}$, $\IHH_n=L^2(0,1)$, $n\in\N\cup\{\infty\}$ endowed with the Lebesgue measure and the obvious choice of continuous sections. We let $\form{h}_n$ be the classical Dirichlet form corresponding to the Laplacian with Neumann boundary conditions for $n\in\N$ and with Dirichlet conditions for $n=\infty$, i.e. $D(\form{h}_n)=W^{1,2}(0,1)$ for $n\in \N$, $D(\form{h}_\infty)=W_0^{1,2}(0,1)$. Clearly, this section defines a continuous field. However, taking the constant function $1$ for all $n\in\N$, we get a section that converges weakly (even strongly) in $\IHH$, but its limit $1\not\in D(\form{h}_\infty)$.  
\end{example}

\section{Examples, convergence of sequences, and nets}
Here we record examples of continuous fields mostly of Hilbert spaces. The important point we make is that all the classes of examples from the literature dealing with notions of strong convergence of operators defined on varying spaces we are aware of can be cast in the framework introduced in the previous sections. 
\begin{example} (Trivial fields)
 Let $X$ be a topological space, $\IVV$ a locally convex space, whose topology is induced by a family $p=(p_i;\; i\in I)$. Then $X\times \IVV:=(\IVV;\; x\in X)$ endowed with
 the constant family $p$ defines a field of locally convex spaces;
$$ \ICC:=\ICC(X\times \IVV):=\Gamma(\IVV)\cap C(X;\IVV)
$$
defines a locally uniformly closed and saturated continuous subspace for $(X\times\IVV,p)$ which makes $(X\times\IVV,\ICC,p)$  a continuous field of locally convex spaces, called the \emph{trivial field} over $X$.
\end{example}

\begin{example} (Ambient space)
Let $X$ be a topological space, let $\KK$ be a Hilbert space and for $x\in X$ let $P(x)$ be an orthogonal projection in $\KK$, whose range $\IHH_x$ is considered as a Hilbert space. If 
$$
X\ni x\longmapsto P(x)\in\ILL(\KK)
$$ 
is strongly continuous, then the induced collection $\IHH$ together with
$$
\ICC_0:=\set{P(\cdot)\zeta}{\zeta\in \KK}
$$
induces a continuous field of Hilbert spaces over $X$ in that $\ICC:=\overline{\ICC_0}$ is a continuous subspace and it holds that 
$$
\ICC=\set{Pf}{f\in C(X;\IKK)} 
$$
and the fibers carry the inner product from $\IKK$. Note that by \cite{dix}, Proposition 11, p. 240, or our generalization thereof in Proposition \ref{prop:subfield}, much less than strong continuity of the orthogonal projections is needed in case the underlying space is paracompact.
\end{example}
We mention that the preceding class of examples allows to treat domain perturbations for regular Dirichlet forms as in \cite{stollman}, where prior results in \cite{rauchtaylor,simon,weidmann} are generalized. With the convergence theorems from Section \ref{sec:self-adjointops}, we can extend the results and at the same time simplify the proofs considerably. We will return to a detailed account in a future paper.

\begin{example}\label{fam} (Continuous families of seminorms on a fixed vector space) We fix a vector space $V$ over $\K$, and consider a topological space $X$ together with family $\net{\|\cdot\|}{x}{X}$ so that for every $v\in V$ the map
$$
X\ni x\longmapsto \| v\|_x\in [0,\infty)
$$ 
is continuous. By dividing out the kernel $\INN_x$, we get a continuous field of normed spaces $(\IEE,\ICC)$ over $X$, where
$$
\IEE_x=V/\INN_x\mbox{  and  }\ICC:=\overline{\set{X\ni x\mapsto [v]_{\IEE_x}\in \IEE_x}{v\in V}}^{\IEE},
$$
the closure in the sense of Lemma \ref{Lemma2.4}. By fiberwise completion, we get a continuous field of Banach spaces $(\tilde{\IEE},\ICC_{\tilde{\IEE}})$ by Proposition \ref{completion2}, where $\ICC_{\tilde{\IEE}}:=\overline{\ICC}^{\tilde{\IEE}}$. If the original seminorms come from semi inner products, the resulting field is a continuous field of Hilbert spaces. In the above situation, we say that \emph{$\ICC_{\tilde{\IEE}}$ is generated by the constant sections $\set{X\ni x\mapsto [v]_{\IEE_x}\in \IEE_x}{v\in V}$}.  
\end{example}

\begin{example} (Tangent bundle of Wasserstein space) We produce a concrete instance of Example \ref{fam}: let $M=(M,g)$ be a smooth connected Riemannian manifold and for $p\in [1,\infty)$ consider the space $\mathscr{P}_p(M)$ of Borel probability measures $\mu$ on $M$ with finite $p$-th moments, that is
$$
\int_M d(x,x_0)^p d\mu(x)<\infty \quad\text{for some/all $x_0\in M$}.
$$
This space is turned into a metric space with respect to the Wasserstein distance
$$
W_p(\mu_1,\mu_2):=\inf\set{\left(\int_{M\times M} d(x,y)^p d\pi(x,y)\right)^{1/p}}{\text{$\pi$ is a coupling of $\mu_1$ and $\mu_2$}}.
$$
With 
$$
\mathfrak{X}^{\mathrm{grad}}_c(M):=\{\nabla f: f\in C^\infty_c(M) \}
$$
the space of smooth compactly supported gradient vector fields on $M$, one defines the tangent space over $\mu \in \mathscr{P}_p(M)$ by
$$
T_\mu \mathscr{P}_p(M):=\overline{\mathfrak{X}^{\mathrm{grad}}_c(M)}^{L^p_\mu(TM)},
$$
where $L^p_\mu(TM)$ denotes the space of $L^p$-vector fields on $M$ with respect $\mu$. As convergence with respect to $W_p$ implies weak convergence (see Theorem 7.12 in \cite{villani}), it follows that the maps
$$
\mu\longmapsto \left\|\nabla f\right\|_{L^p_\mu(TM)}
$$
are continuous. With
$$
T \mathscr{P}_p(M):= \bigsqcup_{\mu\in  \mathscr{P}_p(M)}T_\mu \mathscr{P}_p(M)
$$
we thus obtain a continuous field of Banach spaces 
$$
\big(T \mathscr{P}_p(M),\ICC(T \mathscr{P}_p(M))\big),
$$
where $\ICC(T \mathscr{P}_p(M))$ is generated by the constant sections $\{\mu\mapsto [A]_\mu: A\in\mathfrak{X}^{\mathrm{grad}}_c(M)\}$, where $[A]_\mu$ denotes the equivalence class of the vector field $A$ with respect to $\mu$. For $p=2$ this yields a continuous field of Hilbert spaces, which can be considered as the tangent bundle of $\mathscr{P}_2(M)$, and in this case the space of continuous sections $\ICC(T \mathscr{P}_2(M))$ is given by all $A\in\Gamma(T\mathscr{P}_2(M))$
such that 
\begin{itemize}
\item  for all smooth compactly supported vector fields $B$ on $M$, the map
$$
\mu\longmapsto \left\langle A(\mu),B\right\rangle_{\mu}:=\int_M (A(\mu),B) \;d\mu
$$
is continuous,
\item the map
$$
\mu\longmapsto \left\|A(\mu)\right\|_{\mu}=\sqrt{\left\langle A(\mu),A(\mu)\right\rangle_{\mu}}
$$
is continuous.
\end{itemize}
This characterization follows from Proposition \ref{prop:mwandnormgivesstrong}. Subbundles of this tangent bundle have played a major role in the geometrization of PDEs, ever since the seminal works by Otto \cite{otto} and Benamou/Brenier \cite{brenier}, and also Lott \cite{lott}. If we look at $M=\mathbb{R}^m$ and take $\mu=\delta_0$ the delta measure, we find that $T_\mu \mathscr{P}_p(M)=\mathbb{R}^m$ is finite dimensional, while for smooth $\mu$'s this space is infinite dimensional. Thus, these tangent bundles are a very natural instance of a field or bundle which is not locally trivial in any reasonably continuous sense.
\end{example}

It is our pleasure to thank Melchior Wirth for pointing that the previous example might fall into our framework.

\begin{example}\label{loccomp} ($L^p$-spaces on locally compact spaces; see also \cite{bos}.) Here we start with a locally compact Hausdorff space $T$. The base of the continuous field to be constructed is the space $\IMM^R(T)$ of nonnegative Radon measures on $T$, equipped with the weak-*-topology coming from $C_c(T)$. We can then regard, for any $p\in [1,\infty]$, the continuous field of Banach spaces over $\IMM^R(T)$ given by 
$$
\big((L^p(T,\mu))_{\mu\in \IMM^R(T)},\ICC_{T,p}),
$$
as another concrete instance of Example \ref{fam}, where $\ICC_{T,p}$ is generated by the constant sections $\{\mu\mapsto [\varphi]_\mu\;|\>\varphi\in C_c(T)\}$, with $[\varphi]_\mu$ the equivalence class with respect to the measure $\mu$. For $p=2$ we obtain a continuous field of Hilbert spaces with the following explicit (cf. Proposition \ref{prop:mwandnormgivesstrong}) description of its space of continuous sections: $\ICC_{T,2}$ is given by all sections $f$ of $(L^p(T,\mu))_{\mu\in \IMM^R(T)}$ such that 
\begin{itemize}
\item the map 
$$
\mu\longmapsto \left\langle\varphi, f(\mu)\right\rangle_{\mu}:=\int \varphi f(\mu) \;d\mu 
$$
is continuous for all  $\varphi\in C_c(T)$, 
\item the map
$$
\mu\longmapsto \left\|f(\mu)\right\|_{\mu}=\sqrt{\left\langle f(\mu), f(\mu)\right\rangle_{\mu}} 
$$
is continuous.
\end{itemize}

\end{example}

\begin{example}(Fields of Riemannian manifolds / Ricci flow) Assume $M$ is a smooth connected manifold, $X$ is a metric space and $(g^{(x)})_{x\in X}$, is a continuous family of smooth Riemannian metrics on $M$, in the sense that each $g^{(x)}$ is a smooth Riemannian metric on $M$, and in each chart $U\subset M$, with the corresponding metric coefficients $g^{(x)}_{ij}$, the map
$$
X\times U \longrightarrow \mathbb{R}, \quad (x,y)\longmapsto g^{(x)}_{ij}(y)
$$
is continuous. Let $\mu_x$ be the Riemannian volume measure associated with $g^{(x)}$. \\
In the sequel, we consider 
$$
W:=L^2_\loc(M)_+=\{w\in L^2_\loc(M): w\geq 0\}
$$
as a Fr\'echet space with the topology induced by the seminorms $\|\cdot\|_{L^2(K,\mu)}$, $K\subset M$ compact, where $\mu$ is any fixed Riemannian volume measure on $M$ (the opology is independent of a particular choice). \\
Next we are going to show that this situation induces a continuous family of unbounded operators in the sense of Theorem \ref{Cfield-Mcontinuity}, namely the Friedrichs realizations $T_{\form{t}_{(x,v)}}$ of the underlying Schrödinger operators $-\Delta_{g^{(x)}}+w$ in $L^2(M,\mu_{x})$: let $\form{t}_{(x,w)}$ denote the form which is given by the closure of the form 
\begin{align*}
\form{t}_{(x,w)}[\phi]:=\int_M g^{(x)}_*(d\phi,d\phi) \;d\mu_{x}+\int_M w \phi^2 d\mu_x=-\int_M (\Delta_{g^{(x)}} \phi) \phi \;d\mu_{x}+\int_M w \phi^2 d\mu_x,\quad \phi\in C^\infty_c(M),
\end{align*}
where $g^{(x)}_*$ denotes the co-metric on $T^*M$. Its domain of definition $\IDD(\form{t}_{(x,w)})$ is thus given as the completion of $C^\infty_c(M)$ with respect to the  scalar product
\begin{align*}
\left\langle \phi_1,\phi_2\right\rangle_{\form{t}_{(x,w)}}&:=\left\langle \phi_1,\phi_2\right\rangle_{x} +\form{t}_{(x,w)}[\phi_1,\phi_2].
\end{align*}

As a base of the corresponding field of Hilbert spaces we use $X\times W$ and the field 
$$
\IHH:=(\IHH_{(x,w)};(x,w)\in X\times W),
$$
where 
$$
\IHH_{(x,w)}:=L^2(M,\mu_x).
$$
By pulling back via the map $(x,w)\mapsto \mu_x$, the field $\IHH$ becomes a continuous field of Hilbert spaces with the space of sections
$$
\ICC:=\overline{\set{(x,w)\mapsto f(\mu_x)}{f\in \ICC_{M,2}}},
$$
where $\ICC_{M,2}$ is as in Example \ref{loccomp} and the closure in the sense of Lemma \ref{Lemma2.4}. Note that while the generating sections are all constant with respect to $w\in W$, this is no longer true for elements of $\ICC$.
We will show that 
$$
\form{t}=(\form{t}_{(x,w)};(x,w)\in X\times W)
$$
is Mosco-continuous, using Corollary \ref{cor:moscoseq} (iii). Clearly, (M2') is satisfied:  for all  $\varphi\in C^\infty_c(M)$, the section
$$
f: (x,w)\longmapsto [\varphi]_{\mu_x}\in \IHH_{(x,w)}\subset \IHH
$$ 
is continuous, as is the map
$$
(x,w)\longmapsto \form{t}_{(x,w)}[\varphi],
$$
hence $f\in \ICC(\form{t})$; this gives that $C^\infty_c(M)\subset \set{g(x)}{g\in \ICC(\form{t})}$
for every $x\in X$.
\iffalse
and we define a family of seminorms on $C^\infty_c(M)$ by 
$$
\left\|\varphi\right\|^2_{(x,w)}:=\int_M \varphi^2 \;d\mu_x,\quad (x,w)\in  X\times L^2_\loc(M)_+.
$$

It follows easily from dominated convergence that the map $(x,w)\mapsto \left\|\varphi\right\|_{(x,w)}$ is continuous. We thus get a continuous field of Hilbert spaces 
$$
\big(L^2(M,\mu_{\bullet}),\ICC^{\bullet}_{M}\big)
$$
over $X\times L^2_\loc(M)_+$, whose fiber over $(x,w)$ is given by $L^2(M,\mu_{x})$ and whose space of continuous sections $\ICC^{\bullet}_{M}$ is given by all $f\in\Gamma(L^2(M,\mu_{\bullet}))$ with the following two properties:
\begin{itemize}
\item for all fixed $\varphi\in C^\infty_c(M)$, the map 
$$
(x,w)\longmapsto  \langle f(x,w),\varphi\rangle_{x} := \int_M f(x,w) \varphi \;d\mu_x 
$$
is continuous,
\item the map
$$
(x,w)\longmapsto \left\|f(x,w)\right\|_{x}=\sqrt{\langle f(x,w),f(x,w)\rangle_{x}}
$$
is continuous.
\end{itemize}
Next we are going to show that this situation induces a continuous family of unbounded operators in the sense of Theorem \ref{Cfield-Mcontinuity}, namely the Friedrichs realizations $T_{\form{t}_{(x,v)}}$ of the underlying Schrödinger operators $-\Delta_{g^{(x)}}+w$ in $L^2(M,\mu_{x})$: let $\form{t}_{(x,w)}$ denote the form which is given by the closure of the form 
\begin{align*}
\form{t}_{(x,w)}[\phi]:=\int_M g^{(x)}_*(d\phi,d\phi) \;d\mu_{x}+\int_M w \phi^2 d\mu_x=-\int_M (\Delta_{g^{(x)}} \phi) \phi \;d\mu_{x}+\int_M w \phi^2 d\mu_x,\quad \phi\in C^\infty_c(M),
\end{align*}
where $g^{(x)}_*$ denotes the co-metric on $T^*M$. Its domain of definition $\IDD(\form{t}_{(x,w)})$ is thus given as the completion of $C^\infty_c(M)$ with respect to the  scalar product
\begin{align*}
\left\langle \phi_1,\phi_2\right\rangle_{\form{t}_{(x,w)}}&:=\left\langle \phi_1,\phi_2\right\rangle_{x} +\form{t}_{(x,w)}[\phi_1,\phi_2].
\end{align*}
Actually, one can show (e.g. Theorem XIII.1. in \cite{batu}) 
$$
\IDD(\form{t}_{(x,w)})= W^{1,2}_0(M,g^{(x)})\cap L^2(M,w\;d\mu_x),\quad \text{where $W^{1,2}_0(M,g^{(x)})=\IDD(\form{t}_{(x,w)})\mid_{w=0}$}. 
$$
In particular, one has
$$
\IDD(\form{t}_{x,w})\subset W^{1,2}_\loc(M)\cap L^2_\loc(M;w),
$$
where $L^2_\loc(M;w)$ stands for the weighted $L^2_\loc$-space. In particular, the latter intersection becomes a (metric-independent) Fr\'echet space in the obvious way,
It follows follows from similar arguments as before that 
$$
\IDD(\form{t}_{\bullet}):= \bigsqcup_{(x,w)\in X\times L^2_\loc(M)_+}\IDD(\form{t}_{(x,w)})
$$
becomes a continuous field of Hilbert spaces over $X\times L^2_\loc(M)_+$, whose space of continuous sections given by all $f\in\Gamma(\IDD(\form{t}_{\bullet}))$
such that 
\begin{itemize}
\item for all fixed $\varphi\in C^\infty_c(M)$, the map 
$$
(x,w)\longmapsto (\form{t}_{(x,w)}+1)[f(x,w),\varphi]
$$ 
is continuous,
\item the map 
$$
(x,w)\longmapsto \left\|f(x,w)\right\|_{\form{t}_{(x,w)}}=\sqrt{\langle f(x,w),f(x,w)\rangle_{\form{t}_{(x,w)}}}
$$
is continuous.
\end{itemize}

Any such section is in fact an element of $\ICC^\bullet_{M}$, which is implied by the first statement of Proposition \ref{prop:mwandnormgivesstrong}. In particular, this space coincides with $\mathscr{C}(\form{t}_\bullet)$ defined prior to Theorem \ref{Cfield-Mcontinuity}. 
\else
It remains to show property (M1') from Corollary  \ref{cor:moscoseq} (iii). So assume $f\in \mathscr{C}(\form{t})^\w$ is bounded and let $x_n\to x_\infty$ in $X$ and $w_n\to w_\infty$ in $W$. We set
$$
(\cdot,\cdot)_n:=g^{(x_n)}_*,\quad (\cdot,\cdot)_\infty:=g^{(x_\infty)}_*,\quad \mu_n:=\mu_{x_n},\quad \mu_\infty:=\mu_{x_\infty},
$$
and define for each $n$ a smooth section $T_{n}$ of $T^*M$ with pointwise strictly positive eigenvalues via 
\begin{align}\label{e11}
(\alpha,\beta)_n=(T_{n}\alpha,\beta)_\infty,\quad\text{so that}\quad d\mu_{n}= \det(T_{n})^{\frac{1}{2}}d\mu_{\infty}.
\end{align}
Clearly, $T_n(x)\to 1$ in $\mathrm{End}(T^*_xM)$ for all $x$ (this is metric-independent). Moreover, as all $x_n$ are in a compact subset of $X$ and as $g^{(x)}(y)$ depends continuously on $(x,y)\in X\times M$, for each $K\subset M$ compact one finds a $c>0$ such that for all $n$ and all $y\in K$ one has 
\begin{align}\label{edxaed}
c \; (\cdot,\cdot)_n\leq (\cdot,\cdot)_\infty\leq (1/c)\; (\cdot,\cdot)_n\quad\text{as symmetric bilinear forms in $T_y M$}.
\end{align}
With $f_n:=f(x_n)$ and $f_\infty:=f(x_\infty)$ we have to show that for all $\varphi\in C^\infty_c(M)$ one has 
$$
\int_M f_n \varphi\; d\mu_{n}=\int_M f_n \varphi \det(T_{n})^{\frac{1}{2}}\; d\mu_{\infty}\to \int_M f_\infty \varphi\; d\mu_{\infty}.
$$
It suffices to show
\begin{align}\label{ewqa}
\int_M f_n \varphi \; d\mu_{\infty}\to \int_M f_\infty \varphi\; d\mu_{\infty}.
\end{align}
Indeed, (\ref{ewqa}) is sufficient, because
\begin{align}\label{wdld}
\int_M f_n \varphi \; d\mu_{n}- \int_M f_n \varphi\; d\mu_{\infty}=\int_M f_n \varphi (\det(T_{n})^{\frac{1}{2}}-1)d\mu_{\infty}\to 0,
\end{align}
which in turn with $K:=\mathrm{supp}(\varphi)$ can be seen as follows,
\begin{align*}
& \int_M |f_n \varphi (\det(T_{n})^{\frac{1}{2}}-1)|d\mu_{\infty}\\
&\leq \left\|\varphi\right\|_\infty \left(\int_K |f_n |^2d\mu_{\infty}\right)^{\frac{1}{2}}\left(\int_K |  (\det(T_{n})^{\frac{1}{2}}-1)|^2d\mu_{\infty}\right)^{\frac{1}{2}}\\
&=\left\|\varphi\right\|_\infty\left(\int_K |f_n |^2\det(T_n)^{-\frac{1}{2}}\; d\mu_{n}\right)^{\frac{1}{2}}\left(\int_K |  (\det(T_{n})^{\frac{1}{2}}-1)|^2d\mu_{\infty}\right)^{\frac{1}{2}}\\
&\leq C\left\|\varphi\right\|_\infty\left(\int_K |  (\det(T_{n})^{\frac{1}{2}}-1)|^2d\mu_{\infty}\right)^{\frac{1}{2}},\quad C:=\sup_{n\in\mathbb{N}}\left\|1_K\det(T_n)^{-\frac{1}{2}}\right\|_\infty\left(\int_M f_n^2\; d\mu_{n}\right)^{\frac{1}{2}}<\infty,
 \end{align*}
because $f$ is a bounded section of $\IDD(\form{t})$ and because of (\ref{edxaed}), and finally 
$$
\int_K |  (\det(T_{n})^{\frac{1}{2}}-1)|^2d\mu_{\infty}\to 0,
$$
because of (\ref{edxaed}) and dominated convergence.\\
To prove (\ref{ewqa}), it suffices to prove that 
\begin{align}\label{wqaw}
&\int_M (df_n ,d\varphi)_\infty \; d\mu_{\infty}+\int_M f_n \varphi \; d\mu_{\infty}+\int_M f_n \varphi \, w_\infty\; d\mu_{\infty}\\\nonumber
&\to \int_M (df_\infty ,d\varphi)_\infty \; d\mu_{\infty}+\int_M f_\infty \varphi\; d\mu_{\infty}+\int_M f_\infty \varphi \, w_\infty\; d\mu_{\infty},
\end{align}
because this implies that $f_n$ converges weakly in the Fr\'{e}chet space $W^{1,2}_\loc(M)\cap L^2_\loc(M,w_\infty\,dx)$ and thus weakly in $L^{2}_\loc(M)$ to $f_\infty$, giving (\ref{ewqa}). \\
To prove (\ref{wqaw}), note that as $f$ is a weakly continuous section of $\IDD(\form{t})$, we have 
\begin{align}\label{wqaw2}
&\int_M (df_n ,d\varphi)_n \; d\mu_{n}+\int_M f_n \varphi \; d\mu_{n}+\int_M f_n \varphi\; w_n\; d\mu_{n}\\\nonumber
&\to \int_M (df_\infty ,d\varphi)_\infty \; d\mu_{\infty}+\int_M f_\infty \varphi\; d\mu_{\infty}+\int_M f_\infty \varphi  \;w_\infty\; d\mu_{\infty}.
\end{align}
It is thus enough to show
\begin{align*}
&\int_M f_n \varphi \; d\mu_{\infty}- \int_M f_n \varphi \; d\mu_{n}\to 0,\\
&\int_M (df_n ,d\varphi)_\infty \; d\mu_{\infty}- \int_M (df_n ,d\varphi)_n \; d\mu_{n}\to 0,\\
&\int_M f_n \varphi\; w_\infty\; d\mu_{\infty}-\int_M f_n \varphi\; w_n\; d\mu_{n}\to 0.
\end{align*}
The first difference has already been treated above. The second difference is 
\begin{align*}
=\int_K \left(\big(1-\det(T_n)^{\frac{1}{2}}T_n\big)df_n ,d\varphi\right)_\infty  \; d\mu_{\infty},
\end{align*}
the absolute value of which is 
\begin{align*}
&\leq \left(\sup_K|d\varphi|_\infty\right)\left(\int_K |df_n|^2_\infty   \; d\mu_{\infty}\right)^{\frac{1}{2}}\left(\int_K \|\big(1-\det(T_n)^{\frac{1}{2}}T_n\|^2 \; d\mu_{\infty}\right)^{\frac{1}{2}}\\
&= \left(\sup_K|d\varphi|_\infty\right)\left(\int_K (T_n^{-1}df_n,df_n)_n \det(T_n)^{-\frac{1}{2}}  \; d\mu_{n}\right)^{\frac{1}{2}}\left(\int_K \|\big(1-\det(T_n)^{\frac{1}{2}}T_n\|^2 \; d\mu_{\infty}\right)^{\frac{1}{2}}\\
&\leq C\left(\int_K \|\big(1-\det(T_n)^{\frac{1}{2}}T_n\|^2 \; d\mu_{\infty}\right)^{\frac{1}{2}},
\end{align*}
which goes to $0$, noting that
$$
C:=\left(\sup_K|d\varphi|_\infty\right)\left(\sup_n\sup_K\left\|\det(T_n)^{-\frac{1}{2}}T_n^{-1}\right\|_n\right)\sup_n\left(\int_M |df_n|^2_n \; d\mu_{n}\right)^{\frac{1}{2}}<\infty.
$$

The third difference is
\begin{align*}
= \int_M f_n \varphi\left(w_\infty-w_n\det(T_n)^{\frac{1}{2}}\right)\; d\mu_{\infty},
\end{align*}
the absolute value of which is
\begin{align*} 
&\leq \left(\int_M f_n^2\varphi\; d\mu_\infty\right)^{\frac{1}{2}}  \left(\int_K \left|w_\infty-w_n\det(T_n)^{\frac{1}{2}}\right|^2\; d\mu_\infty\right)^{\frac{1}{2}}\\
&\leq \left\|\varphi\right\|_\infty\left(\int_M f_n^2 \det(T_n)^{-\frac{1}{2}}\; d\mu_n\right)^{\frac{1}{2}}  \left(\int_K \left|w_\infty-w_n\det(T_n)^{\frac{1}{2}}\right|^2\; d\mu_\infty\right)^{\frac{1}{2}}\\
&\leq C \left(\int_K \left|w_\infty-w_n\det(T_n)^{\frac{1}{2}}\right|^2\; d\mu_\infty\right)^{\frac{1}{2}},
\end{align*}
where
$$
C:=\left\|\varphi\right\|_\infty\sup_{n}\left(\left\|1_K\det(T_n)^{-\frac{1}{2}}\right\|_\infty\left(\int_M f_n^2 \; d\mu_n\right)^{\frac{1}{2}} \right)<\infty.
$$
In view of $\det(T_n)^{\frac{1}{2}}\to 1$ pointwise on $K$ and
$$
\sup_{n} \left\|1_K\det(T_n)^{\frac{1}{2}}\right\|_\infty<\infty,
$$
the last term tends to $0$, in view of $w_n\to w_\infty$ in $L^2_\loc(M)$. This completes the proof of assumption (M1') from Corollary \ref{cor:moscoseq} (iii). \\

An important instance of this example is provided by any Ricci flow on $M$, where $X=I$ is a time interval, and $g^{(t)}$ solves
$$
\partial_t g^{(t)}= -2\mathrm{Ric}_{g^{(t)}},\quad t\in I.
$$
\end{example}

\begin{example} (Discrete approximations) We already mentioned the approach by Stummel from the 1970s, see \cite{stummel1,stummel2}: consider normed spaces $E_n$, $n \in\N$, 
and $E$ and a linear mapping 
$$
R:E\longrightarrow \Big(\prod_{n\in\N}E_n\Big)/\sim ,
$$
where $\snet{\zeta}{n}{\N}\sim \snet{\eta}{n}{\N}$ whenever $\| \zeta_n-\eta_n\|\to 0$ as $n\to\infty$. The map $R$ is assumed to satisfy the norm continuity condition
$$
\snet{\zeta}{n}{\N}\in R(\zeta_\infty)\quad\Longrightarrow\quad \lim_{n\to\infty}\|\zeta_n\|=\|\zeta_\infty\| .
$$
We therefore get a continuous field over $X=\N\cup\{\infty\}$ by considering the family of spaces above together with
$$\ICC=\set{\snet{\zeta}{n}{\N\cup\{\infty\}}}{\snet{\zeta}{n}{\N}\in R(\zeta_\infty)}.
$$
\end{example}

\begin{example}(Kuwae-Shioya)
 In \cite{kuwae}, the authors consider nets of operators over varying spaces. With the notation of the present text, let $A$ be a directed set and 
 $$\IHH:=(\IHH_\alpha;\alpha \in A\cup\{\infty\})$$ 
 a field of Hilbert spaces. For any $\alpha\in A$ there is a densely defined linear operator $\Phi_\alpha:C_\infty\to \IHH_\alpha$ with the additional property that $\| \Phi_\alpha u\|\to \| u\|$ for all $u\in C_\infty$. With $\Phi_\infty:=\mathrm{id}:C_\infty\to C_\infty$ we see that 
 $$\ICC_0:=\set{\Phi(\cdot)u}{u\in C_\infty}$$
 is a continuous subspace. Its closure $\ICC$ in the sense of Lemma \ref{Lemma2.4} induces a continuous field by Lemma \ref{completion}. Since the topology on $A$ is discrete, 
 $$
 \ICC^{\mathsf{KS}}:=\set{\zeta\in\Gamma(\IHH)}{\text{there exists $f\in\ICC$ such that $\lim_{\alpha\in A}\|\zeta(\alpha)-f(\alpha)\|=0$}}
 $$
 defines a continuous field as well and the corresponding strong convergence is exactly the one that is used by Kuwae and Shioya. In this sense, our constructions can be considered as generalizations of the results in \cite{kuwae}, as we have mentioned several times already. In particular, our framework applies to the examples treated in the latter paper, which include e.g.
 \begin{itemize}
 \item pointed measured Gromov Hausdorff convergence of metric measure spaces,
 \item Lipschitz convergence of Riemannian manifolds,
 \item Collapsing of warped product metrics,
 \item convergence of discrete graphs.
 \end{itemize}
In addition, \cite{kuwae} has also been specialized to further similar and  other settings, such as smooth settings (e.g. \cite{bei, carron}), singular (metric measure) settings (e.g. \cite{honda,honda2}), discrete settings (e.g. \cite{kasue, lohr}), and infinite dimensional settings (e.g. \cite{dello,grothauswittmann}).
\end{example}

The next example contains a result that is completely new and could equally well be stated as a theorem in its own right. We deal with generalized Schrödinger operators of the form $H+\mu$, with appropriate measures $\mu$, see \cite{abr,blanchard,sv,sturm} and the references there for further literature and background. While our perturbation result deals with operators on a fixed Hilbert space, we use methods from continuous fields:

\begin{example}(Approximation of singularly perturbed Dirichlet forms)
Let $\form{h}$ be a regular Dirichlet form with domain $D(\form{h})\subset L^2(T,m)$ with $m$ being a Radon  measure of full support on the locally compact Hausdorff space $T$. Let $\IMM^R_0(T)$ be the set of positive Radon measures on $X$ which are absolutely continuous with respect to the capacity induced by  $\form{h}$. We endow $\IMM^R_0(T)$ with the weak-*-topology $\sigma(\IMM^R_0(T),C_c(T))$. For $\mu\in\IMM^R_0(T)$ we consider the form $\form{h}+\mu$ in $L^2(T,m)$ given by
\[
D(\form{h}+\mu) := D(\form{h})\cap L^2(T,\mu),\quad  (\form{h}+\mu)[u] = \form{h}[u] + \int_T \tilde{u}^2\,d\mu,
\]
where $\tilde{u}$ is a quasicontinuous version of $u$; then each $\form{h}+\mu$ is a Dirichlet form and $C_c(T)\cap D(\form{h})$ is a core for it, see \cite{sv} for more details; the associated self-adjoint operator is denoted by $H+\mu$.
We will show:

\textbf{Claim:} For every $\mu_\infty\in \IMM^R_0(T)$ that is singular with respect to $m$, 
$$
\IMM^R_0(T)\longrightarrow \ISS(L^2(T,m)), \quad \mu \longmapsto  H+\mu
$$
is continuous at $\mu_\infty$, where $\ISS(L^2(T,m))$ is equipped with the topology of strong resolvent convergence.

That is, we have to show that, for $\mu_\infty$ as above, and $\mu_\alpha\xrightarrow{w^*}\mu_\infty$, $\form{h}+\mu_\alpha\xrightarrow{\Mos}\form{h}+\mu_\infty$. 

(M2*) from Corollary \ref{mosco} follows with constant nets of the form $\eta_\alpha=\phi\in D(\form{h})\cap C_c(T)$.

We now show (M1). To this end we can assume that $u_\alpha\xrightarrow{\w}u_\infty$ and 
$$
\lim_\alpha(\form{h}+\mu_\alpha+1)[u_\alpha]= \tau_\infty <\infty;
$$
we have to deduce that $(\form{h}+\mu_\infty+1)[u_\infty]\le \tau_\infty$.  We can regard $\snet{\tilde{u}}{\alpha}{A}$ as a bounded net in the continuous field $(L^2(T,\mu))_{\mu\in \IMM^R(T)}$ introduced in Example 7.5 above, as well as in the analogously defined  continuous field $(L^2(T,m+\mu))_{\mu\in \IMM^R(T)}$. Using Proposition \ref{prop:bddgivesmw} iteratively, we get a subnet $\varphi:B\to A$ so that
\begin{align*}
    \tilde{u}_{\varphi(\beta)}&\xrightarrow{\m\w}v_\infty\mbox{  in  }\big(L^2(T,\mu)\big)_{\mu\in \IMM^R(T)}\\
    \tilde{u}_{\varphi(\beta)}&\xrightarrow{\m\w}w_\infty\mbox{  in  }\big(L^2(T,m+\mu)\big)_{\mu\in \IMM^R(T)} ;
\end{align*}
since 'constant nets' $\xi_\alpha=\phi\in  C_c(T)$ are strongly convergent in both fields above, we get:
\begin{align*}
    \langle v_\infty, \phi\rangle_{L^2(T,m+\mu_\infty)}&=\lim_\beta\left(\int_T u_{\varphi(\beta)}\phi d m+
    \int_T \tilde{u}_{\varphi(\beta)}\phi d \mu_{\varphi(\beta)}\right)\\
    &=\langle u_\infty, \phi\rangle_{L^2(T,m)} +\langle w_\infty, \phi\rangle_{L^2(T,\mu_\infty)},
\end{align*}
giving 
$$
\langle v_\infty - u_\infty, \phi\rangle_{L^2(T,m)} =\langle w_\infty-v_\infty, \phi\rangle_{L^2(T,\mu_\infty)}\mbox{   for all  }\phi\in C_c(T) .
$$
Since $m$ and $\mu_\infty$ are mutually singular, this gives $v_\infty=u_\infty$ $m$-a.e. and $w_\infty=v_\infty$ $(m+\mu_\infty)$-a.e., so that $w_\infty=\tilde{u}_\infty$ $\mu_\infty$-a.e.

Finally, since $\form{h}$ is closed, we have that
$(\form{h}+1)[u_\infty]\le \liminf_\beta (\form{h}+1)[u_{\varphi(\beta)}]$. Putting everything together gives
\begin{align*}
    (\form{h}+\mu_\infty+1)[u_\infty] &=(\form{h}+1)[u_\infty]+\int_T\tilde{u}_\infty^2 d\mu_\infty\\
    &\le \liminf_\beta \left( (\form{h}+1)[u_{\varphi(\beta)}] + \int_T\tilde{u}_\beta^2 d\mu_\beta\right)=\tau_\infty,
\end{align*}
using the weak lower semicontinuity of the norm,  Corollary \ref{cor:lscnorm}, in the last step.
\end{example}
The convergence result in the previous example can be seen as complementary to Theorem 3 in \cite{listol}, where all measures are assumed to be absolutely continuous w.r.t.~the reference measure, negative parts are allowed which have to satisfy certain uniform relative boundedness conditions, but the convergence is stronger than the weak-*-convergence in the preceding example. For details and references to earlier results see loc.~cit.

\end{document}